\documentclass{article}

\usepackage{tikz}
\usepackage{tikz-cd}
\usepackage{url, hyperref}

\tikzset{->-/.style={decoration={
  markings,
  mark=at position .45 with {\arrow{>}}},postaction={decorate}}}

\usepackage{amsmath}
\usepackage{amssymb}
\usepackage{enumitem}
\usepackage{graphicx}
\usepackage{mathdots}
\usepackage{color}

\def\oM{\overline{\mathcal{M}}}

\def\Z{\mathbb{Z}}

\def\qed{{\hfill $\Diamond$}}
\def\b1{{\bf 1}}

\def\Aut{{\rm Aut}}

\def\E{\mathrm{E}}

\def\L{\mathrm{L}}
\def\V{\mathrm{V}}
\def\H{\mathrm{H}}
\def\g{\mathrm{g}}

\def\P{\mathsf{P}}

\def\com{\mathbb{C}}

\usepackage{theorem}
\newtheorem{definition}{Definition}
\newtheorem{theorem}[definition]{Theorem}

\newtheorem{proposition}[definition]{Proposition}
\newtheorem{corollary}[definition]{Corollary}

\newcommand{\T}{\mathsf{T}}

\newcommand{\Joh}[1]{}
\newcommand{\Rah}[1]{}
\newcommand{\Sam}[1]{}

\begin{document}

\title{The Hodge bundle,  the universal 0-section, and the log Chow ring
of the moduli space of curves}

\author{S. Molcho, R. Pandharipande, J. Schmitt}

\date{September 2022}

\maketitle

\begin{abstract}
  We bound from below
  the complexity of the top Chern class $\lambda_g$ of the Hodge
  bundle in the Chow ring of the moduli space of curves: no
  formulas for $\lambda_g$ in terms of classes of degrees 1 and 2 can exist.
  As a consequence of the Torelli map, the
  0-section over the second Voronoi compactification of
  the moduli of principally polarized abelian varieties
  also can not be expressed in terms  of classes of degree 1 and 2.
  Along the way, we establish new cases of Pixton's conjecture for
  tautological relations.

  In the log Chow ring of the moduli space of curves, however,
  we prove $\lambda_g$  lies
  in the subalgebra generated by logarithmic boundary divisors.
  The proof is effective and uses Pixton's double ramification cycle formula together with a foundational study of the
  tautological ring defined by a normal crossings
  divisor. 
  The results open the door to the search for  simpler
  formulas for $\lambda_g$ on the moduli
  of curves after log blow-ups.
  \end{abstract}

\setcounter{tocdepth}{1}
\tableofcontents

\section{Introduction}
\subsection{The Hodge bundle}
Let $\overline{\mathcal{M}}_g$ be the moduli space of Deligne-Mumford stable curves, and
let
$$\pi: \mathcal{C}_g \rightarrow  \overline{\mathcal{M}}_g$$
be the universal curve with relative dualizing sheaf $\omega_\pi$.
The rank $g$ Hodge bundle
$\mathbb{E}_g$ on
 $\overline{\mathcal{M}}_g$
is
defined by
$$\mathbb{E}_g = \pi_* \omega_\pi\, .$$

The study of the Chern classes of the Hodge bundle goes back at least to
Mumford's Grothendieck-Riemann-Roch calculation  \cite{Mum} in the 1980s.
Starting in the late 1990s, the connection
of the Hodge
bundle to the deformation theory of the moduli space of stable maps has
led to an exploration
of Hodge integrals in various contexts, see \cite{AKMV,ELSV,FPHodge, GrP, LLLZ, LLZ,MOOP, OkP,P}.

The top Chern class{\footnote{All Chow classes are taken here with $\mathbb{Q}$-coefficients.}}
of the Hodge bundle 
$$\lambda_g = c_g(\mathbb{E}_g) \in \mathsf{CH}^g(\overline{\mathcal{M}}_g)$$
plays a special role for several reasons:
\begin{enumerate}
\item[(i)]  Two {\em vanishing properties} hold:
  $$\lambda_g^2=0 \in \mathsf{CH}^{2g}(\overline{\mathcal{M}}_g) \ \ \text{and}\ \ \lambda_g|_{\Delta_0}= 0  \in
  \mathsf{CH}^g(\Delta_0)\,, $$
  where $\Delta_0\subset \overline{\mathcal{M}}_g$ is the
  divisor of curves with a non-separating node.
  The first vanishing follows from the highest graded part of Mumford's relation 
  $$c(\mathbb{E}_g)\cdot c(\mathbb{E}_g^*)= 1\, ,$$
  proven in \cite[equations (5.4), (5.5)]{Mum}. The second follows from the existence of a trivial quotient{\footnote{The
      quotient is defined on the double cover of $\Delta_0$ obtained
      by ordering the branches of the node.}}
  $$\mathbb{E}_g \twoheadrightarrow \mathbb{C}$$
  determined by
  the residue at (a branch of) the node, see
  \cite[Section 0.4]{FPFul}.
\item[(ii)] The class $(-1)^g \lambda_g$ appears
  in the virtual fundamental class of the moduli of 
  {\em contracted maps} in the Gromov-Witten theory of target curves.
  Since the {\em double ramification cycle} in the degree 0 case is defined
    via contracted maps, we have
    $$\mathsf{DR}_{g,(0,\ldots,0)} = (-1)^g\lambda_g \in \mathsf{CH}^g(\overline{\mathcal{M}}_{g,n})\, ,$$
where $\overline{\mathcal{M}}_{g,n}$ is the moduli space of stable pointed
    curves.
    See \cite[Sections 0.5.3 and 3.1]{JPPZ}.

  Another basic consequence
  is the $\lambda_g$-formula \cite{FPLg},
  \begin{equation*}
    \int_{\overline{\mathcal{M}}_{g,n}} \psi_1^{k_1} \ldots \psi_n^{k_n} \lambda_g =
    \binom{2g+n-3}{k_1,\ldots, k_n} \cdot \int_{\overline{\mathcal{M}}_{g,1}} \psi_1^{2g-2} \lambda_g\, ,
    \end{equation*}
  predicted by the Virasoro constraints for degree 0 maps to curves \cite{GP}. Here
    $$\psi_i = c_1(\mathbb{L}_i) \in \mathsf{CH}^1(\overline{\mathcal{M}}_{g,n})$$
    is the Chern class of the cotangent line at the $i^{th}$ point.
    The $\lambda_g$-formula  plays a central role in the
    study of the tautological ring $\mathsf{R}^\star(\mathcal{M}^{\mathsf{ct}}_{g,n})$
    of the moduli space of curves of compact type \cite{P-ct}.

\item[(iii)] Again as an excess class, $(-1)^g\lambda_g$ appears fundamentally in the local Gromov-Witten
theory
of surfaces. For example, the Katz-Klemm-Vafa formula \cite{kkv} proven in \cite{mpt,ptmul} concerns
integrals
$$\int_{[\overline{\mathcal{M}}_g(S,\beta)]^{\mathsf{red}}} (-1)^g \lambda_g$$
against the reduced virtual fundamental class of the moduli space of stable maps
to $K3$ surfaces. For a recent study of the parallel problem for local log Calabi-Yau surfaces (with integrand $(-1)^g \lambda_g$), see \cite{Bouss}.
    
  \item[(iv)] The class $(-1)^g\lambda_g$ arises via the pull-back of the
    universal $0$-section 
    of 
    the moduli space of {\em principally polarized abelian varieties} (PPAVs). Over the moduli space of
    compact type curves, the connection to PPAVs shows a
    third vanishing property:
    $$\lambda_g |_{\mathcal{M}_{g}^{\mathsf{ct}}} = 0\, ,$$
    see \cite{vdG}. We will discuss PPAVs  further
    in Section \ref{zzss} below.
  \end{enumerate}

  Our main results here concern the complexity of the class $\lambda_g$ in the Chow ring.
  For $\overline{\mathcal{M}}_g$, we bound from below
  the complexity of formulas for
  $$\lambda_g\in \mathsf{CH}^\star(\overline{\mathcal{M}}_g)\, .$$
  As a consequence of the connection to the moduli of PPAVs, we also
  bound from below
  the complexity of formulas for the universal $0$-section.

  The log Chow ring of
  $(\overline{\mathcal{M}}_g,\partial \overline{\mathcal{M}}_g)$ is defined
  as a colimit over all iterated blow-ups of boundary strata.
  The usual Chow ring is naturally a subalgebra
  $$\mathsf{CH}^\star(\overline{\mathcal{M}}_g) \subset
  \mathsf{logCH}^\star(\overline{\mathcal{M}}_g,
  \partial \overline{\mathcal{M}}_g)\, .$$
  The main positive result of the paper is the simplicity of $\lambda_g$ in the
  log Chow ring. We prove
  $$\lambda_g \in \mathsf{divlogCH}^\star(\overline{\mathcal{M}}_g,
\partial \overline{\mathcal{M}}_g)\, , $$
where
$$\mathsf{divlogCH}^\star(\overline{\mathcal{M}}_g, \partial \overline{\mathcal{M}}_g)
\subset  \mathsf{logCH}^\star(\overline{\mathcal{M}}_g,
\partial \overline{\mathcal{M}}_g
)$$
  is the subalgebra generated by logarithmic boundary divisors.
While 
$\lambda_g$ in Chow
is complicated, 
$\lambda_g$ in log Chow
is as simple as possible!
 We present several related open questions.
  
\subsection{The 0-section} \label{zzss}

Let $\mathcal{A}_g$ be the moduli space of PPAVs of dimension $g$, and let
$$\pi: \mathcal{X}_g \rightarrow \mathcal{A}_g $$
be the universal abelian variety $\pi$ equipped with a universal $0$-section
$$s: \mathcal{A}_g \rightarrow \mathcal{X}_g\, .$$
The image of the $0$-section determines an algebraic cycle class
$$Z_g \in \mathsf{CH}^g(\mathcal{X}_g)\, .$$
The second Voronoi compactification of $\mathcal{A}_g$  has been given a modular
interpretation by Alekseev:
$$\mathcal{A}_g \subset \overline{\mathcal{A}}^{\mathsf{Alekseev}}_g\, .$$
Olsson \cite{Ols} provided a
modular interpretation for the normalization
$$\overline{\mathcal{A}}^{\mathsf{Olsson}}\rightarrow\, 
\overline{\mathcal{A}}_g^{\mathsf{Alekseev}}\,. $$
Our approach here will be equally valid for both
$\overline{\mathcal{A}}^{\mathsf{Olsson}}$ and
$\overline{\mathcal{A}}_g^{\mathsf{Alekseev}}$.
We will simply denote the compactification by
$$\mathcal{A}_g \subset \overline{\mathcal{A}}_g\, ,$$
where $\overline{\mathcal{A}}_g$ stand for either the space of Alekseev or
the space of Olsson.

The four important properties{\footnote{We follow the notation of  \cite{Ols}.}}
of the compactification
$\overline{\mathcal{A}}_g$ which we will require are:

\begin{enumerate}
\item[$\bullet$] The points of $\overline{\mathcal{A}}_g$ 
  parameterize (before normalization) stable semiabelic pairs which are
  quadruples $(G,P,L,\theta)$ where $G$ is a semiabelian
  variety, $P$ is a projective variety equipped with a $G$-action,
  $L$ is an ample line bundle on $P$, and $\theta\in H^0(P,L)$.
  The data $(G,P,L,\theta)$ satisfy several further conditions, see
  Section 4.2.16 of \cite{Ols}.

\item[$\bullet$] There is a universal semiabelian variety
  $$\overline{\pi}: \overline{\mathcal{X}}_g \rightarrow \overline{\mathcal{A}}_g$$
  with a $0$-section
  $$\overline{s}: \overline{\mathcal{A}}_g \rightarrow \overline{\mathcal{X}}_g\, $$
  corresponding to the semiabelian variety which is
  the first piece of data of a stable semiabelic pair (the rest of
  the pair data will not play a role in our study).
  
\item[$\bullet$] The usual Torelli map $\tau: \mathcal{M}_g
  \rightarrow \mathcal{A}_g$ extends canonically
  $$\overline{\tau}: \overline{\mathcal{M}}_g \rightarrow
  \overline{\mathcal{A}}_g\, ,$$
  see \cite{Alek}.

\item[$\bullet$] The $\overline{\tau}$-pullback to $\overline{\mathcal{M}}_g$
  of $\overline{\mathcal{X}}_g$ is the universal family
  $$\mathsf{Pic}_\epsilon^0 \rightarrow \overline{\mathcal{M}}_g$$
  parameterizing line bundles on the fibers of the universal curve
  $$\epsilon:\mathcal{C}_g \rightarrow \overline{\mathcal{M}}_g$$
  which have degree 0 {\em on every component} of any fiber \cite{Alek}.

\end{enumerate}


The image of the $0$-section $\overline{s}$ determines an
operational Chow class
$$\overline{Z}_g\in \mathsf{CH}^g_{\mathsf{op}}(\overline{\mathcal{X}}_g)$$
since the image is an \'etale local complete intersection
in $\overline{\mathcal{X}}_g$. The class $\overline{Z}_g$ is related to $(-1)^g\lambda_g$ via a pull-back construction.
Let
$$t:\overline{\mathcal{M}}_g \rightarrow \mathsf{Pic}_\epsilon^0 $$
be the $0$-section defined by the trivial line bundle.
By the properties of
$$\overline{\pi}: \overline{\mathcal{X}}_g \rightarrow \overline{\mathcal{A}}_g$$
discussed above,
$$\overline{\tau}^*\overline{s}^*(\overline{Z}_g) =
t^*(t_*[\overline{\mathcal{M}}_g]) 
\, .$$
By the standard analysis of  the vertical tangent bundle of
$\mathsf{Pic}_\epsilon^0$,
$$t^*(t_*[\overline{\mathcal{M}}_g]) = (-1)^g \lambda_g \in
\mathsf{CH}^g(\overline{\mathcal{M}}_g)\, .$$
Indeed, by the excess intersection formula the class $t^*(t_*[\overline{\mathcal{M}}_g])$ equals the top Chern class of the normal bundle of the zero section of $\mathsf{Pic}_\epsilon^0$. Over $[C]\in\overline{\mathcal{M}}_g$, the fiber of the normal bundle
is the first-order deformation space of the trivial line bundle on $C$. 
The deformation space is given by $$H^1(C, \mathcal{O}_C) = H^0(C, \omega_C)^\vee\, ,$$ the fiber of the dual of the Hodge bundle $\mathbb{E}_g^\vee$ with top Chern class $(-1)^g \lambda_g$.
We conclude
\begin{equation}\label{ff66}
  \overline{\tau}^*\overline{s}^*(\overline{Z}_g)  = (-1)^g \lambda_g
\in
\mathsf{CH}^g(\overline{\mathcal{M}}_g)\, .
\end{equation}

\subsection{Complexity of the $0$-section}

The study the $0$-section over $\mathcal{A}_g$
is related to the double ramification cycle (especially over curves of compact type), see  Hain \cite{Hain} and
Grushevsky-Zakharov \cite{GZ1}. A central idea there is to use the beautiful formula
\begin{equation} \label{dd1}
  Z_g = \frac{\Theta^g}{g!}\ \in \mathsf{CH}^g(\mathcal{X}_g)\, , \end{equation}
where $\Theta \in \mathsf{CH}^1(\mathcal{X}_g)$ is the universal symmetric theta divisor trivialized along the $0$-section.
The proof of \eqref{dd1} in Chow uses the Fourier-Mukai transformation and work of Deninger-Murre \cite{DM}, see \cite{BL,V}.
The article \cite{GZ1} provides a more detailed discussion of the history of \eqref{dd1}.

We are interested in the following question:
{\em to what
  extent is an 
  equation of the form of \eqref{dd1}  possible over
  $\overline{\mathcal{A}}_g$?}
A result by Grushevsky and Zakharov along these lines appears in \cite{GZ2}.
As before, let
$$\overline{Z}_g \in \mathsf{CH}^g_{\mathsf{op}}(\overline{\mathcal{X}}_g)$$
be the class of the $0$-section $\overline{s}$.
Grushevsky and Zakharov calculate the restriction  $\overline{Z}_g|_{\mathcal{U}_g}$
of $\overline{Z}_g$ 
over a particular open set{\footnote{$\mathcal{U}_g$ is the locus determined
    by semiabelian varieties of torus rank at most 1.}}
$$\mathcal{A}_g\subset \mathcal{U}_g \subset \overline{\mathcal{A}}_g$$
in terms of $\Theta$, a boundary divisor $D\in\mathsf{CH}^1(\overline{\mathcal{X}}_g|_{\mathcal{U}_g})$, and a class $$\Delta \in \mathsf{CH}^2(\overline{\mathcal{X}}_g|_{\mathcal{U}_g})\, .$$
The result of Grushevsky-Zarkhov shows that 
while the naive extension of \eqref{dd1} does
{\em not} hold over $\mathcal{U}_g$, the class
$\overline{Z}_g|_{\mathcal{U}_g}$
lies in the subalgebra of $\mathsf{CH}^\star(\overline{\mathcal{X}}_g|_{\mathcal{U}_g})$
generated by classes of degrees 1 and 2.
The formula of \cite{GZ2} is
a useful extension of \eqref{dd1}.

The divisor classes $\mathsf{CH}^1_{\mathsf{op}}(\overline{\mathcal{X}}_g)$
generate a subalgebra
$$\mathsf{divCH}^\star_{\mathsf{op}}(\overline{\mathcal{X}}_g)
\subset \mathsf{CH}^\star_{\mathsf{op}}(\overline{\mathcal{X}}_g)\, .$$ 
The first  bound from below of the complexity of the class of
the $0$-section is the following result.

\begin{theorem} \label{ooo} For all $g\geq 3$, we have
$\overline{Z}_g\notin
\mathsf{divCH}^\star_{\mathsf{op}}(\overline{\mathcal{X}}_g)$.

\end{theorem}

 As a consequence, no divisor formula extending
\eqref{dd1} is possible for $\overline{\mathcal{A}}_g$. 
Though not stated,
    the analysis of \cite{GZ2}  over $\mathcal{U}_g$ 
    can be used to show  $\overline{Z}_g|_{\mathcal{U}_g}$ is {\em not}
    in the subalgebra of $\mathsf{CH}^\star(\overline{\mathcal{X}}_g|_{\mathcal{U}_g})$
    generated by classes of degree 1. Theorem 1 can therefore
also be obtained from \cite{GZ2}.{\footnote{We thank S. Grushevsky for correspondence
    about \cite{GZ2}.}}

In fact, we can go further. Let 
$$\mathsf{CH}^\star_{\leq k}(\overline{\mathcal{X}}_g)
\subset \mathsf{CH}^*_{\mathsf{op}}(\overline{\mathcal{X}}_g)\, $$ 
be the subalgebra generated by all elements of degree
at most $k$, so
$$\mathsf{divCH}^\star_{\mathsf{op}}(\overline{\mathcal{X}}_g)
= \mathsf{CH}^\star_{\leq 1}(\overline{\mathcal{X}}_g)\, .$$

\begin{theorem} \label{ttt}  For all $g\geq 7$, we have
$\overline{Z}_g\notin
\mathsf{CH}^\star_{\leq 2}(\overline{\mathcal{X}}_g)$.

\end{theorem}

By Theorem \ref{ttt}, the Grushevsky-Zakharov formula for $\overline{Z}_g
|_{{\mathcal{U}}_g}$ will
require corrections by higher degree classes
when extended over $\overline{\mathcal{A}}_g$. We propose the
following conjecture about the complexity of the class $\overline{Z}_g$.

\vspace{10pt}
\noindent{\bf Conjecture A.} {\em  No extension of \eqref{dd1} over
  $\overline{\mathcal{A}}_g$ for all $g$
 can be written in terms of classes of uniformly bounded degree.}
 \vspace{10pt}
 
 
 The pull-back relation \eqref{ff66} relates the
complexity of 
the class $$\lambda_g\in \mathsf{CH}^\star(\overline
{\mathcal{M}}_g)\, $$ 
to the complexity of 
 $\overline{Z}_g\in
\mathsf{CH}^\star_{\mathsf{op}}(\overline{\mathcal{X}}_g)$.
Theorems \ref{ooo} and \ref{ttt} will be
immediate consequence of 
parallel{\footnote{In fact, we will prove
     in Section \ref{tttw} stronger results in cohomology instead of Chow.}}
complexity bounds for $\lambda_g$.

\subsection{Complexity of $\lambda_g$}

The divisor classes $\mathsf{CH}^1(\overline{\mathcal{M}}_g)$
generate a subalgebra
$$\mathsf{divCH}^\star(\overline{\mathcal{M}}_g)
\subset \mathsf{CH}^\star(\overline{\mathcal{M}}_g)\, .$$ 
The first  bound from below of the complexity of $\lambda_g$
is the following result.

\begin{theorem} \label{rrr} For all $g\geq 3$, we have
$\lambda_g\notin
\mathsf{divCH}^\star(\overline{\mathcal{M}}_g)$.

\end{theorem}

Via the pull-back relation \eqref{ff66}, Theorem \ref{rrr} immediately implies Theorem \ref{ooo}. The
proof of Theorem \ref{rrr}, presented in Section \ref{tttw}, starts with explicit calculations
in the tautological ring in genus $3$ and $4$ using the Sage package
{\em admcycles} \cite{adm}. A boundary restriction argument is then used 
to inductively control all higher genera.

For the analogue of Theorem \ref{ttt}, let
$$\mathsf{CH}^\star_{\leq k}(\overline{\mathcal{M}}_g)
\subset \mathsf{CH}^*(\overline{\mathcal{M}}_g)\, $$ 
be the subalgebra generated by all elements of degree
at most $k$. A similar strategy
(with a much more complicated initial calculation in genus 5)
yields the following result which implies Theorem \ref{ttt}.

\begin{theorem} \label{fff}  For all $g\geq 7$, we have
$\lambda_g\notin
\mathsf{CH}^\star_{\leq 2}(\overline{\mathcal{M}}_g)$.

\end{theorem}

The proofs of Theorems \ref{rrr}  and
\ref{fff} require new cases of Pixton's conjecture about
the ideal of relations in the tautological ring
$$\mathsf{R}^\star(\overline{\mathcal{M}}_{g,n}) \subset
\mathsf{CH}^\star(\overline{\mathcal{M}}_{g,n})\, .$$

\begin{proposition} Pixton's relations generate all relations
  among tautological classes in \label{Pixxx} $\mathsf{R}^4(\overline{\mathcal{M}}_{4,1})$
  and $\mathsf{R}^5(\overline{\mathcal{M}}_{5,1})$.
  \end{proposition}

  While the above arguments become harder to pursue in general for
  $\mathsf{CH}^\star_{\leq k}(\overline{\mathcal{M}}_g)$, we
  expect the following to hold.
  
\vspace{10pt}
\noindent{\bf Conjecture B.} {\em  For fixed $k$,
  $\lambda_g \in \mathsf{CH}^\star_{\leq k}(\overline{\mathcal{M}}_g)$
holds only for finitely many $g$.}
 \vspace{10pt}
 
\noindent Of course, Conjecture B implies Conjecture A.

\subsection{Log Chow}
Theorems \ref{ooo}-\ref{fff}
about the classes  $\overline{Z}_g$ and $\lambda_g$ 
are in a sense negative results since formula types are excluded.
Our main positive result about $\lambda_g$ concerns the larger log Chow ring
$$\mathsf{CH}^\star(\overline{\mathcal{M}}_g) \subset
\mathsf{logCH}^\star(\overline{\mathcal{M}}_g,
\partial\overline{\mathcal{M}}_g)\, . $$
The log Chow ring 
and the subalgebra $$\mathsf{divlogCH}^\star
(\overline{\mathcal{M}}_g,
\partial\overline{\mathcal{M}}_g)$$
generated by logarithmic boundary divisors
are defined carefully in
Section \ref{lcr}. Our perspective,
using limits over log blow-ups, requires
the least background in log geometry.
A more intrinsic
 approach to the definitions can be found
 in \cite{Barrott}.

\begin{theorem}\label{xxx}
  For all $g\geq 2$, we have $\lambda_g \in \mathsf{divlogCH}^\star(\overline{\mathcal{M}}_g,
  \partial\overline{\mathcal{M}}_g)
  $.
\end{theorem}

Our proof of Theorem \ref{xxx} 
is constructive: we start with Pixton's formula for the double ramification
  cycle for constant maps \cite{JPPZ} and show each term lies in
  $\mathsf{divlogCH}^\star(\overline{\mathcal{M}}_g)$. In principle, bounds
  for the necessary log blow-ups are possible to obtain from the proof, but these will certainly
  not be optimal. Finding a minimal (or efficient) sequence of log-blows of $(\overline{\mathcal{M}}_g,
  \partial\overline{\mathcal{M}}_g)$ after which $\lambda_g$ lies in the
  subalgebra of logarithmic boundary divisors is an interesting question.

  A crucial part of the proof of Theorem \ref{xxx}
  is the study in Section \ref{dslc} 
  of the logarithmic tautological
  ring, 
  $$\mathsf{R}^\star(X,D) \subset \mathsf{CH}^\star(X)\, ,$$
  defined by a normal crossings divisor 
  $D\subset X$
  in a nonsingular variety $X$.
 Tautological classes are 
 defined here using the Chern roots of the
 normal bundle of logarithmic strata
 $S\subset X$. The precise definitions
 are given in Section \ref{deflc}.
  
  We prove
  three main structural results about logarithmic
  tautological classes:
  \begin{enumerate}
      \item [(i)] $\mathsf{R}^\star(X,D) \subset 
      \mathsf{divlogCH}^\star(X,D)$,
      \item[(ii)]
      pull-backs of tautological
      classes under log blow-ups
      are tautological,
      \item[(iii)] push-forwards
      of tautological classes
      under log blow-ups are
      tautological.
  \end{enumerate}
Our first proof of (i) is presented in Section \ref{lcr2} via an explicit analysis of {\em explosions}: sequences
of blow-ups associated to logarithmic strata of $X$.
A second approach to (i-iii), via the
geometry of the 
  Artin fan of $(X,D)$, is given in 
  Section \ref{craf}. The Artin
  fan perspective, advocated{\footnote{See Ranganathan's 
  April 2020 lecture 
  {\em Gromov-Witten theory and logarithmic intersections}
  in the {\em Algebraic Geometry
  and Moduli Zoominar} at ETH Z\"urich. A foundational
  development will appear in \cite{MoRa}.}}
  by D. Ranganathan, 
  is theoretically  more
  flexible.
  
  After 
  Pixton's formula for the double ramification
  cycle for constant maps
  is shown to lie in
  $\mathsf{R}^\star(\overline{\mathcal{M}}_g,
  \partial\overline{\mathcal{M}}_g)$,
  property (i) implies
  Theorem \ref{xxx}.
  Since Pixton's formula and
  the proof of (i) are both effective, 
  divisor expressions for $\lambda_g$
  are possible to compute {in principle}.
  The result reveals
   the essential simplicity of $\lambda_g$ and
  opens the door to the search for a simpler formula in divisors. 
  
  The proof of Theorem \ref{xxx} yields a refined
  result: only logarithmic boundary divisors over
  $$\Delta_0\subset \overline{\mathcal{M}}_g$$
  are needed to generate $\lambda_g$.
  The parallel result is also
  true for pointed curves:
  $$\lambda_g \in \mathsf{divlogCH}^\star(\overline{\mathcal{M}}_{g,n},\Delta_0)$$
  for $2g-2+n>0$.

  We have seen  that $(-1)^g \lambda_g$ is a special case of the double ramification cycle.
  The general 
  double ramification cycle
  $${\mathsf{DR}}_{g,A} \in \mathsf{CH}^g(\overline{\mathcal{M}}_{g,n})$$
  is defined with respect to a vector of integers 
  $A=(a_1,\ldots,a_n)$  satisfying
  $$\sum_{i=1}^n a_i =0\, .$$
  In \cite[Appendix A]{HPS}, the double ramification cycle was lifted to log Chow{\footnote{The paper \cite{HPS} is primarily formulated in the language of the related bChow ring, which we discuss below and treat in detail in Section \ref{bch}.}},
\begin{equation}\label{ldr}
  \widetilde{\mathsf{DR}}_{g,A} \in \mathsf{logCH}^g(\overline{\mathcal{M}}_{g,n})\, .
  \end{equation}
Motivated by Theorem \ref{xxx}, we conjecture{\footnote{In developments
after the paper was completed, Conjecture C has been proven, see
Section \ref{mgdrc} for a discussion.}} 
a uniform divisorial
property of the lifted double ramification cycle \eqref{ldr}.

\vspace{10pt}
\noindent{\bf Conjecture C.} {\em  For all $g$ and $A$, we have 
  $\widetilde{\mathsf{DR}}_{g,A} \in \underline{\mathsf{div}}\mathsf{logCH}^\star(\overline{\mathcal{M}}_{g,n})$
 where
$$ \underline{\mathsf{div}}\mathsf{logCH}^\star(\overline{\mathcal{M}}_{g,n}) \subset
\mathsf{logCH}^\star(\overline{\mathcal{M}}_{g,n})$$
is the subalgebra generated by logarithmic 
boundary divisors together with the cotangent
line classes $\psi_1,\ldots, \psi_n$.}
   \vspace{10pt}

 Finally, we return to the $\Theta$-formula \eqref{dd1} for $Z_g$.
 {\em Is an extension of the $\Theta$-formula possible over
 $\overline{\mathcal{M}}_g$
 in $\mathsf{logCH}^\star (\overline{\mathcal{M}}_g)$?}
 More specifically, can we find
 $$\mathsf{T}\in \mathsf{logCH}^1(\overline{\mathcal{M}}_g)$$
 which satisfies the following two properties?
 \begin{enumerate}
 \item[(i)] The restriction of $\mathsf{T}$ over the moduli
   of curves $\mathcal{M}^{\mathsf{ct}}_g$ of compact type is $0$.
 \item[(ii)] 
   $(-1)^g\lambda_g = \frac{\mathsf{T}^g}{g!} \in
   \mathsf{logCH}^g(\overline{\mathcal{M}}_g)\, .$
   \end{enumerate}
 Property (i) is imposed since
   $$\Theta|_{Z_g} = 0 \in \mathsf{CH}^1(Z_g)$$
   by the trivialization condition for $\Theta$.
   Unfortunately, the answer is {\em no} even for genus $2$.

   \begin{proposition}\label{noT}
     There does not exist a class $\mathsf{T}\in \mathsf{logCH}^1(\overline{\mathcal{M}}_2)$ satisfying the restriction property (i) and
     $$(-1)^2\lambda_2 = \frac{\mathsf{T}^2}{2!} \in
   \mathsf{logCH}^2(\overline{\mathcal{M}}_2)\, .$$
\end{proposition}

The $\Theta$-formula for $(-1)^g\lambda_g$
can {\em not} be
extended in a straightforward way in $\mathsf{CH}^g(\overline{\mathcal{M}}_g)$
or $\mathsf{logCH}^g(\overline{\mathcal{M}}_g)$. However,
$$\lambda_g\in \mathsf{logCH}^g(\overline{\mathcal{M}}_g)$$
{\em is} a degree $g$ polynomial in the logarithmic boundary divisors over $\Delta_0\subset
\overline{\mathcal{M}}_g$.

\vspace{10pt}
\noindent{\bf Question D.} {\em Find a polynomial formula in
  logarithmic boundary divisors  for $\lambda_g$ in log Chow
  (without using Pixton's formula).}
 \vspace{10pt}

The larger bChow ring  of
$\overline{\mathcal{M}}_g$
is defined as a
limit over {\em all}
blow-ups:
$$\mathsf{CH}^\star(\overline{\mathcal{M}}_g) \subset
\mathsf{logCH}^\star(\overline{\mathcal{M}}_g,
\partial\overline{\mathcal{M}}_g)
\subset \mathsf{bCH}^\star(\overline{\mathcal{M}}_g)
\, .$$
The bChow ring is by far the largest
of the three Chow constructions. 
In Section \ref{bch}, we show
the main questions of the paper become trivial in bChow. 
In fact, for every nonsingular variety $X$, we have
$$\mathsf{divbCH}^\star(X) = \mathsf{bCH}^\star(X)\, .$$
The logarithmic geometry of 
$\overline{\mathcal{M}}_g$ is therefore 
the natural place to
study Question D for $\lambda_g$.

\subsection{Acknowledgments}
D. Holmes,
D. Ranganathan, and J. Wise have suggested that the $\Theta$-formula \eqref{dd1}
should extend over the moduli of curves
in some form in log geometry (based on their understanding
of the logarithmic Picard stack \cite{MoWi}). Our initial motivation here was to
study geometric
obstructions to such an extension. While the simplest form
is excluded, Theorem \ref{xxx} supports the idea of the existence of
some perturbed extension of \eqref{dd1} in log Chow.
Our development of the logarithmic
tautological ring of $(X,D)$ emerged
from the proof of Theorem \ref{xxx}.
We are very grateful to
Holmes, Ranganathan, and Wise for extensive discussions of
these topics.

We have also benefitted from related conversations with Y. Bae,
C. Faber, T. Graber,  S. Grushevsky, M. Olsson,
A. Pixton, R. Vakil, and
D. Zakharov. The
results of the paper were presented in the {\em Algebraic
Geometry Seminar} at Stanford in the fall of 2020 (with a lively and very helpful
discussion afterwards).

Finally, we want to thank the referee for a very careful reading and many suggestions that helped to improve the text.

S.M. was supported by ERC-2017-AdG-786580-MACI.
R.P. was supported by SNF-200020-182181,  ERC-2017-AdG-786580-MACI, and SwissMAP.  J.S. was supported by the SNF Early Postdoc Mobility grant 184245 and thanks
the Max Planck Institute for Mathematics in Bonn for its hospitality. 
This project has received funding
from the European Research Council (ERC) under the European Union Horizon 2020 research and innovation program (grant agreement No 786580).

\section{\texorpdfstring{$\lambda_g$}{lambda\_g} in the Chow ring} \label{tttw}
\subsection{Proof of Theorem \ref{rrr}} \label{tttw1}
{
Recall that the tautological rings $(R^\star(\overline{\mathcal{M}}_{g,n}))_{g,n}$ are defined as the smallest system of $\mathbb{Q}$-subalgebras with unit of the Chow rings $(\mathsf{CH}^\star(\overline{\mathcal{M}}_{g,n}))_{g,n}$ closed under pushforwards by gluing and forgetful maps (see \cite{FPFul, Pcalc} for more details).}
 The tautological
subring $\mathsf{RH}^\star(\overline{\mathcal{M}}_{g,n})$
is defined as the image of the cycle map
$$\mathsf{R}^\star(\overline{\mathcal{M}}_{g,n}) \twoheadrightarrow
\mathsf{RH}^\star(\overline{\mathcal{M}}_{g,n})
\subset \mathsf{H}^{2\star}(\overline{\mathcal{M}}_{g,n})\, .$$
We will use the complex degree grading for $\mathsf{RH}^\star$ and
the real degree grading (as usual) for $\mathsf{H}^\star$.
Let
$$  \mathsf{divRH}^\star(\overline{\mathcal{M}}_{g,n}) \subset 
\mathsf{RH}^{\star}(\overline{\mathcal{M}}_{g,n})\ \ \ {\text{and}} \ \ \
\mathsf{divH}^\star(\overline{\mathcal{M}}_{g,n}) \subset 
\mathsf{H}^{2\star}(\overline{\mathcal{M}}_{g,n})
$$
be the subrings generated respectively by
$\mathsf{RH}^1(\overline{\mathcal{M}}_{g,n})$ and
$\mathsf{H}^2(\overline{\mathcal{M}}_{g,n})$. Since
$$\mathsf{RH}^1(\overline{\mathcal{M}}_{g,n}) =
\mathsf{H}^2(\overline{\mathcal{M}}_{g,n})\, ,$$
by \cite[Theorem 2.2]{arbarellocornalba}
we have
\begin{equation}\label{nn34}
\mathsf{divRH}^\star(\overline{\mathcal{M}}_{g,n})  = 
\mathsf{divH}^{2\star}(\overline{\mathcal{M}}_{g,n})\, .
\end{equation}
We will use the complex degree grading for both 
$\mathsf{divRH}^\star$ and $\mathsf{divH}^\star$.
Since
$$ \mathsf{CH}^1(\overline{\mathcal{M}}_{g,n}) \stackrel{\sim}{=} \mathsf{H}^2(\overline{\mathcal{M}}_{g,n})$$
via the cycle class map, we obtain a surjection
$$ \mathsf{divCH}^\star(\overline{\mathcal{M}}_{g,n})\twoheadrightarrow
\mathsf{divH}^{\star}(\overline{\mathcal{M}}_{g,n})\subset 
\mathsf{H}^{2\star}(\overline{\mathcal{M}}_{g,n})\, .$$
The following stronger result implies Theorem \ref{rrr}.

\vspace{10pt}
\noindent {\bf Theorem 3/Cohomology.}\,  For all $g\geq 3$, we have
$\lambda_g\notin
\mathsf{divH}^\star(\overline{\mathcal{M}}_g)$ .
\vspace{10pt}

\noindent {\em Proof.} For $g=3$, we have complete control of the
tautological rings in Chow and cohomology
since the intersection pairing to 
$\mathsf{R}_0(\overline{\mathcal{M}}_g) \stackrel{\sim}{=}
\mathbb{Q}$ is nondegenerate
for
tautological classes (see \cite{Faber}). In particular,
$$\mathsf{R}^\star(\overline{\mathcal{M}}_3) \stackrel{\sim}{=}
\mathsf{RH}^\star(\overline{\mathcal{M}}_3)\, .$$
In degree 3,
$$\mathsf{divRH}^3(\overline{\mathcal{M}}_3) \subset
\mathsf{RH}^3(\overline{\mathcal{M}}_3)$$
is a 9-dimensional subspace of a 10-dimensional space.
Explicit calculations with the Sage program {\em admcycles} \cite{adm}
show $\lambda_3 \notin \mathsf{divRH}^3(\overline{\mathcal{M}}_3)$.
We conclude $\lambda_3 \notin \mathsf{divH}^\star(\overline{\mathcal{M}}_3)$
by \eqref{nn34}.

Adding one marked point, we can consider the case of $\overline{\mathcal{M}}_{3,1}$. Again it is known that all (even) cohomology classes on $\overline{\mathcal{M}}_{3,1}$ are tautological (see \cite[Section 5.1]{SvZ}). Thus again by Poincar\'e duality, the intersection pairing on $\mathsf{RH}^*(\overline{\mathcal{M}}_{3,1})$  is perfect and hence we can completely identify these groups in terms of generators and relations. One finds that 
$$\mathsf{divRH}^3(\overline{\mathcal{M}}_{3,1}) \subset
\mathsf{RH}^3(\overline{\mathcal{M}}_{3,1})$$
is a 28-dimensional subspace of a 29-dimensional space.
But
remarkably, a calculation by {\em admcycles} shows
\begin{equation*}
  \lambda_3 \in \mathsf{divRH}^3(\overline{\mathcal{M}}_{3,1})\ !
\end{equation*}
The containment appears miraculous. Is there a geometric explanation?

  The tautological ring $\mathsf{RH}^*(\overline{\mathcal{M}}_{4,1})$
  is also completely under control in codimension 4:
$$\mathsf{divRH}^4(\overline{\mathcal{M}}_{4,1}) \subset
\mathsf{RH}^4(\overline{\mathcal{M}}_{4,1})$$
is a 103-dimensional subspace of a 191-dimensional space. An
{\em admcycles} calculation shows
\begin{equation} \label{zz1}
  \lambda_4 \notin \mathsf{divRH}^4(\overline{\mathcal{M}}_{4,1})\, .
  \end{equation}
  The result \eqref{zz1} implies
  $\lambda_4 \notin \mathsf{divRH}^4(\overline{\mathcal{M}}_{4})$ by a pull-back argument and
  $$\lambda_4 \notin \mathsf{divH}^\star(\overline{\mathcal{M}}_{4})$$
  since divisor classes are tautological.

For $g\geq 5$, a  boundary restriction argument is pursued.
Suppose, for contradiction,
\begin{equation}\label{ee1}
  \lambda_g \in \mathsf{divH}^g(\overline{\mathcal{M}}_g)\, .
\end{equation}
Then, by pull-back, we have
\begin{equation}\label{yyz1}
  \lambda_g \in \mathsf{divH}^g(\overline{\mathcal{M}}_{g,1})\, .
\end{equation}

Consider the standard boundary inclusion
$$\delta: \overline{\mathcal{M}}_{g-1,1} \times \overline{\mathcal{M}}_{1,2}
\rightarrow
\overline{\mathcal{M}}_{g,1}\, .$$
As usual, we have
\begin{equation} \label{eqn:deltapullbacksplitting}
\delta^*(\lambda_g) = \lambda_{g-1} \otimes \lambda_1 \, .
\end{equation}
Then \eqref{yyz1} implies
\begin{equation}\label{tt5}
\lambda_{g-1} \otimes \lambda_1 \in
\mathsf{divH}^g( \overline{\mathcal{M}}_{g-1,1} \times \overline{\mathcal{M}}_{1,2} )\, .
\end{equation}
Since $\mathsf{H}^1(\overline{\mathcal{M}}_{g-1,1})$ and  $\mathsf{H}^1(\overline{\mathcal{M}}_{1,2})$
both vanish,
$$
\mathsf{divH}^\star( \overline{\mathcal{M}}_{g-1,1} \times \overline{\mathcal{M}}_{1,2} ) =
\mathsf{divH}^\star( \overline{\mathcal{M}}_{g-1,1})
\otimes \mathsf{divH}^*(\overline{\mathcal{M}}_{1,2} )\, .$$
We therefore can write $\mathsf{divH}^g( \overline{\mathcal{M}}_{g-1,1} \times \overline{\mathcal{M}}_{1,2} )$ as
\begin{eqnarray} \label{fred}
  & &\ \ \ \mathsf{divH}^{g}( \overline{\mathcal{M}}_{g-1,1})
\otimes
\mathsf{divH}^{0}
         (\overline{\mathcal{M}}_{1,2} )
  \\ \nonumber
&\oplus&
\mathsf{divH}^{g-1}( \overline{\mathcal{M}}_{g-1,1}) \otimes
\mathsf{divH}^{1}
         (\overline{\mathcal{M}}_{1,2} )\\ \nonumber
  &\oplus&
\mathsf{divH}^{g-2}( \overline{\mathcal{M}}_{g-1,1}) \otimes
\mathsf{divH}^{2}
           (\overline{\mathcal{M}}_{1,2} ) \, . 
\end{eqnarray}
Since by \eqref{eqn:deltapullbacksplitting} the degree of $\delta^*(\lambda_g)$ splits as $(g-1)+1$ on the two factors, we conclude 
\begin{align*}
 \lambda_{g-1}\otimes \lambda_1 &\in \mathsf{divH}^{g-1}( \overline{\mathcal{M}}_{g-1,1}) \otimes
\mathsf{divH}^{1}
         (\overline{\mathcal{M}}_{1,2})\\ \implies  \lambda_{g-1} &\in \mathsf{divH}^{g-1}(\overline{\mathcal{M}}_{g-1,1})\, ,
\end{align*}
using that $\lambda_1 \neq 0 \in \mathsf{divH}^{1}(\overline{\mathcal{M}}_{1,2})$.
By descending induction, we contradict \eqref{zz1}.
Therefore \eqref{yyz1} and hence also \eqref{ee1} must be false.
\qed

\subsection{With marked points}
\label{wmp}

The proof of  Theorem 3 in cohomology shows
\begin{equation}\label{fpp3}
  \lambda_{g} \notin \mathsf{divH}^{g}(\overline{\mathcal{M}}_{g,1})\, 
\end{equation}
for $g\geq 4$. By using \eqref{fpp3} as a starting point, we can study
$$
\lambda_{g} \in \mathsf{divH}^{g}(\overline{\mathcal{M}}_{g,n})
$$
for $g\geq 4$ and $n\geq 2$ using the boundary restrictions
$$\widehat{\delta}: \overline{\mathcal{M}}_{g,n-1}
\times \overline{\mathcal{M}}_{0,3} \rightarrow
\overline{\mathcal{M}}_{g,n}\, .$$
The argument used in the proof then easily yields the following statement with
markings.

\vspace{10pt}

\noindent {\bf Theorem 3/Markings.}\,  For all $g\geq 4$ and $n\geq 0$, we have
$$\lambda_g\notin
\mathsf{divH}^\star(\overline{\mathcal{M}}_{g,n})\, .$$

\subsection{Proof of Theorem \ref{fff}}

Define the subalgebra of tautological classes
$$\mathsf{RH}^\star_{\leq k}(\overline{\mathcal{M}}_{g,n}) \subset 
\mathsf{RH}^\star(\overline{\mathcal{M}}_{g,n})$$
generated by 
classes of complex degrees less than or equal to $k$. Since
all divisors are tautological,
$$\mathsf{divRH}^\star(\overline{\mathcal{M}}_{g,n}) =
\mathsf{RH}^\star_{\leq 1}(\overline{\mathcal{M}}_{g,n})\, .$$

The arguments in Sections \ref{tttw1} and \ref{wmp}
naturally generalize to address the following
question: {\em when is}
$$\lambda_{g-r} \in \mathsf{RH}^{g-r}_{\leq k}(\overline{\mathcal{M}}_{g,n})\, ?$$

A crucial case of the question (from the point of view of
boundary restriction arguments) is for $n=1$. Let $\mathsf{Q}_g(r,k)$
be the statement
$$ \lambda_{g-r} \notin  \mathsf{RH}^{g-r}_{\leq k}(\overline{\mathcal{M}}_{g,1})$$
which may be true or false.

For example, $\mathsf{Q}_g(r,g-r)$ is false essentially by definition. In fact,
$$\mathsf{Q}_g(s,g-r) \ \  {\text{is false for all $s\geq r$}} $$
for the same reason.
In fact, depending on the parity of $g-r$ it is also false for $s$ slightly below $r$:
$$\mathsf{Q}_g(r-1,g-r)\ \ {\text{is false whenever $g-r$ is odd}}.$$
To see this, note that the even Chern character $\mathrm{ch}_{g-(r-1)}(\mathbb{E}_g)$ vanishes by \cite[Corollary (5.3)]{Mum}. Expressing it in terms of Chern classes $\lambda_i = c_i(\mathbb{E}_g)$ using Newton's identities, we have
\[
0 = \mathrm{ch}_{g-(r-1)}(\mathbb{E}_g) = \frac{(-1)^{g-r+1}}{(g-r+1)!} \lambda_{g-r+1} + \left(\text{polynomial in }\lambda_1, \ldots, \lambda_{g-r} \right)\,.
\]
This proves that $\lambda_{g-r+1}$ can be written in terms of tautological classes of degrees $1, \ldots, g-r$, showing $\mathsf{Q}_g(r-1,g-r)$ to be false.

The boundary arguments used in Sections \ref{tttw1} and \ref{wmp}
yield the following two
    results.

\begin{proposition}\label{cc3}  If $\mathsf{Q}_g(r,k)$ is true,
then  $\mathsf{Q}_{g+1}(r,k)$ and $\mathsf{Q}_{g+1}(r+1,k)$ are true.
\end{proposition}

\begin{proposition}  \label{cc4} If $\mathsf{Q}_g(r,k)$ is true,
then
$$ \lambda_{g-r} \notin  \mathsf{RH}^{g-r}_{\leq k}(\overline{\mathcal{M}}_{g,n})$$
for all $n\geq 0$.
\end{proposition}

Since the $k=1$ case has already been analyzed,
we consider now $k= 2$. The first relevant {\em admcycles}
calculation is 
$$\lambda_3 \notin  \mathsf{RH}^{3}_{\leq 2}(\overline{\mathcal{M}}_{4,1})\, ,$$
so $\mathsf{Q}_4(1,2)$ is true. The corresponding
subspace here is of dimension 91 inside
a 93 dimensional space. As a consequence of Propositions \ref{cc3} and
\ref{cc4} , we
obtain the following result.

\begin{proposition} For all $g\geq 4$ and $n\geq 0$, we have
$$ \lambda_{g-1} \notin  \mathsf{RH}^{g-1}_{\leq 2}(\overline{\mathcal{M}}_{g,n})\, .$$
\end{proposition}

A much more complicated {\em admcycles}
calculation shows
$$\lambda_5 \notin  \mathsf{RH}^{5}_{\leq 2}(\overline{\mathcal{M}}_{5,1})\, ,$$
so $\mathsf{Q}_5(0,2)$ is true. The corresponding
subspace here is of dimension 1314 inside
a 1371 dimensional space.
As a consequence of Propositions \ref{cc3} and \ref{cc4}, we
find
\begin{equation}\label{ncnc}
  \lambda_{g} \notin  \mathsf{RH}^{g}_{\leq 2}(\overline{\mathcal{M}}_{g,n})
  \end{equation}
for all $g\geq 5$ and $n\geq 0$.
For $g\geq 7$, the  equality 
$$\mathsf{RH}^2(\overline{\mathcal{M}}_g) = \mathsf{H}^4(\overline{\mathcal{M}}_g)\, $$
is shown by 
combining results
of Edidin  \cite{Edidin} and Boldsen \cite{Boldsen}.
We provide a summary of the argument in Appendix \ref{cod2}.
For $g\geq 7$, the cycle map 
$$\mathsf{CH}_{\leq 2}^\star(\overline{\mathcal{M}}_g)
\rightarrow \mathsf{H}^{2\star}(\overline{\mathcal{M}}_g)$$
therefore factors through
$\mathsf{RH}_{\leq 2}^\star(\overline{\mathcal{M}}_g)$.
Then, the non-containment \eqref{ncnc} completes
the
proof of Theorem \ref{fff}. \qed

\subsection{Cases of Pixton's conjecture (Proposition \ref{Pixxx})}
For the proofs of Theorem 3 and 4, dimensions and bases of
the following graded parts of tautological rings are required:
\begin{eqnarray*}
\mathsf{RH}^4(\overline{\mathcal{M}}_{4,1})\, , & & 
\text{dim}_{\mathbb{Q}} = 191\, , \\
\mathsf{RH}^5(\overline{\mathcal{M}}_{5,1})\,  , & & 
\text{dim}_{\mathbb{Q}} = 1314\, .
\end{eqnarray*}
These cases are possible to analyze (via {\em admcycles})
since the dual pairings are found
to have kernels exactly 
spanned by Pixton's relations.
A discussion of the 
{\em admcycles} calculation
is presented in Appendix B.
 
Pixton has conjectured  that his relations always provide
all tautological relations. 
Dual pairings are known to be insufficient to prove
Pixton's conjecture in all cases, see \cite{Pcalc, PPZ}
for a more complete discussion.

\section{The  log Chow ring}
\label{lcr}

\subsection{Definitions}
\label{lcrd}

Let $(X,D)$ be a nonsingular variety{\footnote{For a nonsingular
    Deligne-Mumford stack X and a normal crossings divisor $D\subset X$,
    the definitions are the same.}} $X$
with a normal crossings divisor
$$ D = D_1\cup \ldots \cup D_\ell\subset X$$
with $\ell$ irreducible components. The divisor $D\subset X$ is
called the {\em logarithmic boundary}.
An {\em open stratum}
$$S\subset X$$ is an irreducible quasiprojective subvariety satisfying two
properties:
\begin{enumerate}
  \item[(i)]
$S$ is \'etale locally the
transverse intersections of the branches of the $D_i$ which meet $S$.
\item[(ii)] $S$ is maximal with respect to (i).
\end{enumerate}
The set $U=X\setminus D$ is an open stratum. Every
open stratum is nonsingular.
A {\em closed stratum} is the closure of an open stratum.

If all $D_i$ are nonsingular and all
intersections
$$D_{i_1} \cap \ldots \cap D_{i_k}$$ are
irreducible and nonempty,
then there are exactly $2^\ell$ open strata.

Our main  interest will be the 
case $(\overline{\mathcal{M}}_{g,n}, \partial \overline{\mathcal{M}}_{g,n})$
where the normal crossings divisors  have self-intersections.
The open strata defined above for 
 $(\overline{\mathcal{M}}_{g,n}, \partial \overline{\mathcal{M}}_{g,n})$
 are the same as the usual open strata of the
 moduli space of stable curves.

An open stratum $S\subset X$ is {\em simple} if the closure
$$\overline{S} \subset X$$
is
nonsingular.
 A {\em simple blow-up} of $(X,D)$ is a blow-up of $X$ along
 the
 closure  $\overline{S}\subset X$
 of a simple stratum. Let
\begin{equation}\label{nn399}
  \widetilde{X} \rightarrow X
  \end{equation}
 be a simple blow-up along $\overline{S}$. Let
 $$ \widetilde{D}= \widetilde{D}_1 \cup \ldots \widetilde{D}_\ell \cup E \subset
 \widetilde{X}$$
 be the union of the strict transforms $\widetilde{D}_i$ of $D_i$
 along with the exceptional divisor $E$ of the blow-up \eqref{nn399}.
 Then, $(\widetilde{X},\widetilde{D})$ is also a nonsingular
 variety with a normal crossings divisor.
 An {\em iterated blow-up}
 $$(\widehat{X},\widehat{D}) \rightarrow (X,D)$$
 is a finite sequence of simple blow-ups of varieties with
 normal crossings divisors.\footnote{An iterated blow-up is a special type of
 log blow-up. Since we are taking
 a limit, we do not have to
 consider all log blow-ups.}

 The log Chow group of $(X,D)$ is
 defined
 as a colimit over all iterated blow-ups,
\[
\mathsf{logCH}^*(X,D) = \varinjlim_{Y \in \mathsf{logB}(X,D)} \mathsf{CH}^*(Y)\, .
\]
Here, $\mathsf{logB}(X,D)$ is the category of iterated blow-ups of $(X,D)$:
objects in $\mathsf{logB}(X,D)$ are  iterated blow-ups of $(X,D)$
and  morphisms in
$\mathsf{logB}(X,D)$
are iterated blow-ups.

Since $(X,D)$ is the trivial iterated blow-up of itself, there is
canonical algebra homomorphism
$$\mathsf{CH}^\star(X) \rightarrow \mathsf{logCH}^\star(X,D)\, $$
which is injective (since
an inverse map of $\mathbb{Q}$-vectors spaces is obtained by proper
push-forward). We therefore view $\mathsf{CH}^\star(X)$
as a subalgebra of $\mathsf{logCH}^\star(X,D)$. Every
Chow class on $X$ canonically determines a log Chow class for $(X,D)$.

\subsection{Calculation in genus 2}

We will prove Proposition \ref{noT}:
{\em there does not exist a class $\mathsf{T}\in \mathsf{logCH}^1(\overline{\mathcal{M}}_2)$ satisfying
  $$
\mathsf{T}|_{\mathcal{M}^{\mathsf{ct}}_2}=0\ \ \text{and} \ \
  \lambda_2 = \frac{\mathsf{T}^2}{2!} \in
   \mathsf{logCH}^2(\overline{\mathcal{M}}_2)\, .$$}

\noindent{\em Proof.}
Denote by $\pi_* : \mathsf{logCH}^\star(\overline{\mathcal{M}}_{2}) \to
\mathsf{CH}^\star(\overline{\mathcal{M}}_{2})$ the push-forward from log Chow to ordinary Chow.  We will prove a stronger claim: {\em there does not exist
a class
$\mathsf{T}\in \mathsf{logCH}^1(\overline{\mathcal{M}}_{2})$
satisfying
\begin{equation}\label{k99k}
\mathsf{T}|_{\mathcal{M}^{\mathsf{ct}}_2}=0\ \ \text{and} \ \ 
\pi_* \left( \lambda_2 - \frac{\mathsf{T}^2}{2!} \right) = 0
\, \in \mathsf{CH}^2(\overline{\mathcal{M}}_{2})\,.
\end{equation}}

Denote by $U_2 \subseteq \overline{\mathcal{M}}_{2}$ the open subset obtained by removing all closed strata of codimension at least $3$. By the excision exact sequence of Chow groups, we have
\[\mathsf{CH}^2(U_2) \cong \mathsf{CH}^2(\overline{\mathcal{M}}_{2})\]
and thus we can verify the stronger claim by working over $U_2$.

The open set $U_2$ has open strata of codimension 1 and 2.
Since blow-ups along codimension 1 strata do not change $U_2$, 
the only simple blow-ups
$$U_2'\rightarrow U_2$$
are along
codimension 2 open strata (all of which are special in $U_2$). Since the codimension 2 open strata of $U_2$
do not intersect (nor self-intersect), we obtain a $\mathbb{P}^1$-bundle
as an exceptional divisor which contains $0$ and $\infty$ sections{\footnote{Depending upon monodromy, there are either two distinct sections or a
    single double-section{, i.e. a closed subset whose map to the base is finite of degree two}.}}
which are codimension 2 strata of $U'_2$.
The iterated blow-ups
$$\widehat{U}_2 \rightarrow U_2$$ are then simply towers of blow-ups
of these codimension 2 toric strata in successive exceptional divisors.

Assume $\mathsf{T} \in
\mathsf{logCH}^1(U_{2})$ satisfies the conditions \eqref{k99k}.
Since $\mathsf{T}$ restricts to zero over the compact type locus,
$\mathsf{T}$ can be represented as
\[
\mathsf{T} \in \mathsf{CH}^1(\widehat{U}_{2})
\]
on an iterated
blow-up $$\widehat{U}_{2}\rightarrow U_{2}$$ with all blow-up centers living over strata in the complement of the compact type locus. 

\begin{figure}[htb]
    \centering
\begin{tikzpicture}[thick,main node/.style={circle,draw,font=\Large,scale=1.3}]
\node[main node] (A) {0};
\draw (A) to [out=240, in=200,looseness=10] (A);
\draw (A) to [out=160, in=120,looseness=10] (A);
\node (B) at (0,-1.7) {\Large $B$};
\end{tikzpicture}    
\begin{tikzpicture}[thick,main node/.style={circle,draw,font=\Large,scale=1.3}]
\node[main node] (A) {0};
\node[main node] (B) at (2,0) {1};
\draw (A) to [out=200, in=160,looseness=10] (A);
\draw (A) -- (B);
\node (C) at (1,-1.7) {\Large $C$};
\end{tikzpicture}    
    \caption{The stable graphs associated to the codimension $2$ boundary strata $B,C$ contained in $U_2$}
    \label{fig:strataMbar2}
\end{figure}
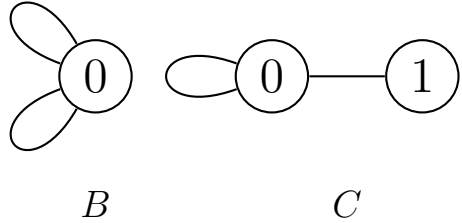

There is a single codimension 1 stratum $\Delta_0\subset U_2$
and  two codimension 2 strata $B,C \subset U_2$ contained in the
complement of the compact type locus (see Figure \ref{fig:strataMbar2}).

Denote by $E_B^1, \ldots, E_B^\ell$ and $E_C^1, \ldots, E_C^m$ the exceptional divisors of blow-ups with centers lying over $B,C$. Then $\mathsf{T}$ has a representation{\footnote{Here, $[\Delta_0]$ is defined via pull-back (not strict
    transformation).}}
\begin{equation*}
    \mathsf{T} = a \cdot [\Delta_0] + \sum_{i=1}^\ell b_i [E_B^i] + \sum_{j=1}^m c_j [E_C^j]\,.
\end{equation*}
After taking the square and pushing forward, we claim
\begin{equation} \label{eqn:Tsquarepushforward}
    \pi_* \left( \mathsf{T}^2 \right) = x \cdot [\Delta_0]^2 + y \cdot [B] + z \cdot [C] \,,
\end{equation}
with $x, y, z \in \mathbb{Q}$ satisfying $$x = a^2 \geq 0\ \ {\text{and}}
\ \ z \leq 0\, .$$
The claim follows from the following observations:
\begin{itemize}
    \item In $\mathsf{T}^2$, all mixed terms $[\Delta_0] \cdot [E_B^i]$ and $[\Delta_0] \cdot [E_C^j]$ vanish after pushforward to $U_{2}$, since
    \[\pi_*([\Delta_0] \cdot [E_B^i]) = [\Delta_0] \cdot \pi_* [E_B^i] =[\Delta_0] \cdot 0 = 0\,. \]
  \item Similarly, since $B \cap C = \emptyset$ in
    $U_2$ (as we have
    removed the codimension $3$ stratum of $\overline{\mathcal{M}}_2$),
    we have $[E_B^i] \cdot [E_C^j] = 0$.
    \item Denote by $\textbf{M} \in \mathsf{Mat}_{\mathbb{Q}, m \times m}$ 
    the matrix defined by
    \[\pi_*\left([E_C^{j_1}] \cdot [E_C^{j_2}] \right) =  \textbf{M}_{j_1, j_2} [C]\,. 
    \]
    A basic fact is that $\textbf{M}$ is negative definite (see  \cite[Section 1]{MumSurface}). Therefore,  for $\textbf{b} = (b_i)_{i=1}^\ell$,
    we have
    \[\pi_* \left( \sum_{j=1}^m b_j [E_C^j] \right)^2 = \underbrace{( \textbf{b}^\top \textbf{M} \textbf{b} )}_{= z \leq 0} [C]\,.\]
    \item The pushforward
    \[\pi_* \left( \sum_{i=1}^\ell b_i [E_B^i] \right)^2\]
    is supported on $B$ and thus is a multiple $y \cdot [B]$ of the fundamental class of $B$.
  \end{itemize}
  
  After substituting  \eqref{eqn:Tsquarepushforward} in the
  second condition of  \eqref{k99k}, we conclude the
  existence of $x, y, z \in \mathbb{Q}$ with $x \geq 0$ and $z \leq 0$
  satisfying
\begin{equation} \label{eqn:lambda2equation}
    x \cdot [\Delta_0]^2 + y \cdot [B] + z \cdot [C] = 2 \lambda_2 \in \mathsf{CH}^2(U_2)\, .
\end{equation}
Using {\em admcycles} (see Appendix \ref{app:lambda2calculation}),
we can explicitly identify all classes in \eqref{eqn:lambda2equation} in
$$\mathsf{CH}^2(U_2) \cong \mathbb{Q}^2\, .$$
The corresponding affine linear equation has the solution space
\[
x = z - \frac{1}{120},\,\ \ y = -\frac{5}{24}\cdot z + \frac{11}{2880}\,.
\]
But for $z \leq 0$, we have
$$z - 1/120 < 0\, ,$$
which contradicts the assumption $x \geq 0$.
Therefore, there can not exist a class $$\mathsf{T} \in
\mathsf{logCH}^1(U_{2})$$ satisfying conditions \eqref{k99k}.
\qed

\section{Relationship with logarithmic geometry}

\subsection{Overview} 
The definitions of Section \ref{lcr} are natural from 
the perspective of logarithmic geometry. The choice of the divisor $D$ on $X$ can be seen as the choice of a log structure on $X$. We briefly recall the relevant definitions and constructions
of logarithmic geometry.

\subsection{Definitions}
A log structure on a scheme $X$ is a sheaf of monoids $M_X$ on the \'etale site of $X$ together with a 
homomorphism{\footnote{$\mathcal{O}_X$ here is sheaf of
monoids under multiplication. In particular, the map $\exp$ transfers the addition in $M_X$ to the usual multiplication of functions in $\mathcal{O}_X$, motivating its name.}}
$$
\textup{exp}:M_X \rightarrow \mathcal{O}_X
$$

\noindent which induces an isomorphism $\textup{exp}^{-1}(\mathcal{O}_X^*) \cong \mathcal{O}_X^*$ on units. 
\begin{enumerate}
    \item[$\bullet$] 
Morphisms of log schemes $(X,M_X) \rightarrow (Y,M_Y)$ are morphisms of schemes $$f:X \to Y$$ together with homomorphisms of sheaves of monoids $f^{-1}M_Y \to M_X$ which are compatible with the structure map $f^{-1}\mathcal{O}_Y \to \mathcal{O}_X$ in the obvious sense. 
\item[$\bullet$]
Log structures can be pulled back. Given a morphism
of schemes $$f: X \to Y\,$$ and a log structure $M_Y$ on $Y$, there is an induced log structure $f^*M_Y$ on $X$, generated by $f^{-1}M_Y$ and the units $\mathcal{O}_X^*$. 
\end{enumerate}
The basics of log schemes can be found in Kato's original article on the subject [Ka].

The category of log schemes is, in practice, too large for geometric study. It is therefore common to work in smaller categories by requiring additional properties
to hold. For our purposes, we will work only with in the category of  fine and saturated
log schemes, usually
termed  {\em f.s. log schemes}. The prototype of such a log scheme is $$A_P = \textup{Spec} (k[P])\,,$$ the spectrum of the algebra generated by a {\em fine and saturated monoid} $P$:
a finitely generated monoid $P$ which injects into its Grothendieck group $P^\textup{gp}$ and which is saturated there, 
$$nx \in P\ \ \text{for}  \ \ n\in \mathbb{N}\, , \
x \in P^\textup{gp} \ \ \ 
\implies \ \ \
x \in P\, .$$ The sheaf $M_{A_P}$ here is the subsheaf of $\mathcal{O}_{A_P}$ generated by $P$ and the units of $\mathcal{O}_{A_P}$.  

All of the log schemes which
arise for us will be 
comparable to $A_P$ on the level of log structures. 
More precisely we require  our log schemes $X$ to admit 
the following local charts: for
each $x \in X$, there must be an \'etale neighborhood $$i:U \to X\, ,$$
an f.s. monoid $P$, and a map $g: U \to A_P$ such that $$i^*M_X = g^*M_{A_P}\,.$$
Since we are always working with f.s. log schemes, the chart $P$ at $x$ can in fact always be chosen to be isomorphic to the characteristic monoid\footnote{Since we are working with \'etale sheaves, the stalk is computed in the \'etale topology; $\overline{x}$ denotes the \'etale stalk.}
$$\overline{M}_{X,\overline{x}}= M_{X,\overline{x}}/\mathcal{O}_{X,\overline{x}}^*$$ at $x$.  

\subsection{Normal crossings pairs}

Let us now return to the situation of interest for the paper: a pair $(X,D)$ of a nonsingular scheme (or Deligne-Mumford stack) with a normal crossings divisor
$D\subset X$. 
The pair $(X,D)$ determines a sheaf $M_X$ on the \'etale site of $X$ by setting 
$$
M_X(p:U \to X) = \{f \in \mathcal{O}_U: f \textup{ is a unit on } p^{-1}(X-D)\}
$$
for each \'etale map $p:U \rightarrow X$. The sheaf of units $\mathcal{O}_X^*$ is a subsheaf of $M_X$. We write $$\overline{M}_X = M_X/\mathcal{O}_X^*$$
for the \emph{characteristic monoid} of $X$.  Normal crossings pairs $(X,D)$, with the log structure described above, are precisely the log schemes which are log smooth over the base field $\textup{Spec } k$ with trivial log structure.  

 When the irreducible components of $D$ do not have self intersections,
the log structure $M_X$ of $(X,D)$ can be defined on the Zariski topology of $X$. 
The result is a technically simpler theory. The pair $(X,D)$ is then
called a {\em toroidal embedding (without self intersection)} in \cite{KKMSD}. However, for a general pair $(X,D)$, $M_X$ can only be defined on the \'etale site of $X$. The general \'etale case differs from the Zariski case in two key aspects: the irreducible components of $D$ can self-intersect, and the characteristic monoid $\overline{M}_X$, while locally constant on a stratum, can globally acquire monodromy. 

The characteristic monoid $\overline{M}_X$ is a constructible sheaf on $X$. The connected components of the loci on which $\overline{M}_X$ is locally constant define
a stratification of $X$, which is precisely the stratification of Section \ref{lcrd}. 
Indeed, for a geometric point $x \in X$, 
$$
\overline{M}_{X,\overline{x}}  = \mathbb{N}^r
$$
where $r$ is the number of branches (in the \'etale topology) of $D$ that contain $x$.

A combinatorial space can be built
from the information contained in $\overline{M}_{X}$. There are two basic approaches. The first, which is more geometric and more evidently combinatorial, is to build the {\em cone complex} $C(X,D)$ of $(X,D)$. We briefly outline the construction (details can be found in
\cite{CCUW} and \cite{ACMUW}). 

We begin with the case where $M_X$ is defined Zariski locally on $X$ (when the irreducible components of $D$ do not have self-intersections). Then, $C(X,D)$ is a rational polyhedral cone complex, see \cite{KKMSD}: 
\begin{enumerate}
\item[$\bullet$] For each point $x \in X$,
the characteristic monoid $\overline{M}_{X,\overline{x}}$ determines a rational polyhedral cone 
$$
\sigma_{X,x} = \textup{Hom}_{\textup{Monoids}}(\overline{M}_{X,\overline{x}}, \mathbb{R}_{\ge 0})
$$
together with an integral structure 
$$
N_{X,x} = \textup{Hom}(\overline{M}_{X,\overline{x}}^\textup{gp},\mathbb{Z})
$$
\item[$\bullet$]
When $x$ belongs to a stratum $S \subset X$ and $y$ belongs to the closure $\overline{S}\subset X$, there are canonical inclusions
$$\sigma_{X,x} \subset \sigma_{X,y}\, , \ \ \ N_{X,x} \subset N_{X,y}\, .$$ 
\item[$\bullet$] We glue the cones $\sigma_{X,x}$ together with their integral structures to form the complex $$C(X,D) = \varinjlim_{x \in X} (\sigma_{X,x}, \sigma_{X,x} \cap N_{X,x})\, .$$ 
\item[$\bullet$] More effectively, instead of working with all points $x \in X$, we can take the finite set $\{x_S \}$ of the generic points of the strata of $(X,D)$. Then, $$C(X,D) = \varinjlim_{x_S} (\sigma_{X,x_S}, \sigma_{X,x_S} \cap N_{X,x_S})\, . $$ In other words, $C(X,D)$ is the dual intersection complex of $(X,D)$. 
\end{enumerate}

When $M_X$ is defined only on the \'etale site, we build the cone complex $C(X,D)$ by descent. 

\begin{enumerate}
\item[$\bullet$] We find an \'etale (but not necessarily proper), strict ($f^*M_X = M_Y)$ cover $f: Y \to X$ which is {\em as fine as possible} (called atomic or small in the literature): the log structure on $Y$ is defined on the Zariski site of $Y$, and each connected component of $Y$ has a unique closed stratum. Taking a further such cover $V$ of the fiber product $Y \times_X Y$ if necessary, we find a groupoid presentation 
$$V \rightrightarrows Y \rightarrow X \, .$$
\item[$\bullet$] We define
$$C(X,D) = \varinjlim [C(V) \rightrightarrows C(Y)] $$
in the category of stacks (with respect to the topology generated by face inclusions) over cone complexes. The construction is carried out in detail in \cite{CCUW}, where it is also shown that it is independent of the choice of groupoid presentation. 
\end{enumerate}

Moreover, $C(X,D)$
is a complex of cones, but 
no longer a rational polyhedral cone complex. For each point $x \in X$, there is a canonical map $$\sigma_{X,x} \rightarrow C(X,D)\, ,$$ 
but the map may no longer be injective. As the \'etale local branches of the divisor $D$ may be connected globally on $X$, the faces of the cones $\sigma_{X,x}$ may be glued to each other in $C(X,D)$, and they may naturally acquire automorphisms coming from the monodromy of the branches of $D$.  

\subsection{Artin fans} \label{afan}
An equivalent combinatorial space is the Artin fan
$\mathcal{A}_X$ of $(X,D)$. The
Artin fan
is defined
by gluing, instead of the dual cones $\sigma_{X,x}$ of $\overline{M}_{X,\overline{x}}$, the quotient stacks 
$$
\mathcal{A}_{\overline{M}_{X,\overline{x}}}=\left[\, \textup{Spec} (k[\overline{M}_{X,\overline{x}}])\,/\, \textup{Spec} (k[\overline{M}_{X,\overline{x}}^\textup{gp}])\,  \right]\, .
$$
 The gluing is exactly the same as for $C(X,D)$ as
 explained above. When $M_X$ is defined on the Zariski site of $X$, 
$$
\mathcal{A}_X = \varinjlim_{x \in X} \mathcal{A}_{\overline{M}_{X,\overline{x}}} = \varinjlim_{x_S}  \mathcal{A}_{\overline{M}_{X,\overline{x_S}}}\, ,
$$
and when $M_X$ is defined only on the \'etale site of $X$,
 $$\mathcal{A}_X = \varinjlim [\mathcal{A}_{V} \rightrightarrows \mathcal{A}_Y]\, ,$$ 
for an atomic presentation $\varinjlim [V \rightrightarrows Y] = X$ as before.

The Artin fan $\mathcal{A}_X$ captures exactly the same combinatorial information as the cone complex $C(X,D)$, but is geometrically less intuitive. Nevertheless, the
Artin fan has the advantage of
coming with a \emph{smooth} morphism of stacks
$$
\alpha:X \rightarrow \mathcal{A}_X\, .
$$

\subsection{Logarithmic modifications}
The cone complex $C(X,D)$ encodes an important operation:  \emph{logarithmic modification} of $X$. 
Logarithmic modifications correspond to subdivisions of $C(X,D)$. A subdivision of $C(X,D)$ is, by definition, a compatible subdivision of all the cones $\sigma_{X,x}$ compatible with the gluing relations. Each cone in the subdivision $\sigma_{X,x}' \rightarrow \sigma_{X,x}$ determines dually a map $\overline{M}_{X,\overline{x}} \rightarrow \overline{M}_{X,\overline{x}}'$, and so a map
$$
\left[\, \textup{Spec} (k[\overline{M}'_{X,\overline{x}})]\, /\, \textup{Spec} (k[\overline{M}_{X,\overline{x}}^\textup{gp}])\, \right] \rightarrow \left[\, \textup{Spec} (k[\overline{M}_{X,\overline{x}}])\, /\, \textup{Spec} (k[\overline{M}_{X,\overline{x}}^\textup{gp}])\,  \right]\, .
$$
The compatibility of the subdivisions with respect to the gluing relations in $C(X,D)$ implies that these maps glue to a \emph{proper}
and {\em birational} representable map 
$$
\mathcal{A}_X' \to \mathcal{A}_X\, .
$$
Then, we define
$$ X' =X \times_{\mathcal{A}_X} \mathcal{A}_X'\rightarrow  X$$
which is proper,
birational, and representable over $X$.
Moreover, $X'$ has an induced log structure, and there is a map
$$\mathcal{A}_{X}' \to \mathcal{A}_{X'}\, $$
which is proper, Deligne-Mumford type, \'etale and bijective. 

The map $\mathcal{A}_{X}' \to \mathcal{A}_{X'}$ -- called the relative Artin fan of $X'$ over $X$ in the litterature -- is not necessarily representable, as the various monodromy groups of the strata of $\mathcal{A}_X$ may act non-faithfully on the strata of $\mathcal{A}_{X}'$, whereas the monodromy groups of the strata of $X'$ act faithfully on $\mathcal{A}_{X'}$ by definition. The strata of $\mathcal{A}_{X}'$ become this way trivial gerbes over the strata of $\mathcal{A}_{X'}$. In a sense, $\mathcal{A}_{X'}$ can be considered as a relative coarse moduli space for $\mathcal{A}_{X}'$\footnote{In fact, this can be made precise: $\mathcal{A}_{X'}$ is the relative coarse moduli space of $\mathcal{A}_{X}'$ with respect to the map $\mathcal{A}_{X}' \to \mathbf{Log}$ to the stack parametrizing log structures.}

Geometrically, subdivisions come in three levels of generality:
\begin{enumerate}
\item[$\bullet$] General subdivisions simply produce proper birational maps $X' \to X$, which are isomorphisms over $X - D$. Such maps are called \emph{logarithmic modifications}
\item[$\bullet$] Log blow-ups are a special kind of subdivision. They are the subdivisions of $C(X,D)$ into the domains of linearity of a piecewise linear function on $C(X,D)$, and they correspond to a sheaf of monomial ideals, $$I \subset M_X\, .$$ The map $X' \to X$ is then projective and is the normalization of the blow-up of $X$ along the sheaf of ideals $\textup{exp}(I) \subset \mathcal{O}_X$. 
\item[$\bullet$] Star subdivisions along simple strata $S$ correspond to the
most basic logarithmic
modifications.
The strata of $X$ are, by construction, in bijection with the cones of $C(X,D)$. We obtain a subdivision by subdividing $\sigma_{X,x_S}$ along its barycenter {(see \cite[Definition 3.3.13]{CLS})}. A simple blow-up along $\overline{S}$ corresponds precisely to the star subdivision of the cone $\sigma_{X,x_S}$. Further applications of the star subdivision operation are discussed in section \ref{expbar}. 
\end{enumerate}

Although star subdivisions are the simplest and most basic subdivisions, we need not consider 
more general subdivisions for our purposes. We are only concerned with statements that are valid over some arbitrarily fine subdivision, and the star subdivisions along simple strata are cofinal in this setting: 
for each subdivision $$C(X,D)' \to C(X,D)\, $$ there is a further subdivision $C(X,D)'' \to C(X,D)'$ such that the
composition $C(X,D)'' \to C(X,D)$ is the composition of star subdivisions along simple strata (see \cite[Chapter 1.7]{Oda}). So the reader can restrict attention to simple blow-ups without any loss of generality. 


We define a category $\mathsf{logM}(X,D)$ whose objects are log modifications $$X' \to X$$ 
obtained via subdivisions of $C(X,D)$. There is a unique morphism $X'' \to X'$ if and only if $X''$ is a log modification of $X'$. Following \cite{Barrott}, we then define
$$
\mathsf{logCH}^\star(X,D) = \varinjlim_{X' \in \mathsf{logM}(X)}\mathsf{CH}^\star(X') \, .$$
As simple blowups are cofinal among log modifications, we have, equivalently, 
$$\mathsf{logCH}^\star(X,D) = \varinjlim_{X'\in \mathsf{logB}(X,D)}\mathsf{CH}^\star(X')\, $$
as defined in Section \ref{lcrd}.

\section{The divisor subalgebra of log Chow}
\label{dslc}

\subsection{Definitions} \label{deflc}
Let $(X,D)$ be a nonsingular variety $X$ with a normal crossings divisor
$$ D = D_1\cup \ldots \cup D_\ell\subset X$$
with $\ell$ irreducible components.
Let
$$\mathsf{divlogCH}^\star(X,D) \subset \mathsf{logCH}^\star(X,D)$$
be the subalgebra
generated by the classes of all the components 
of the associated normal crossings divisors of all
iterated blow-ups of $X$.

 Let $S\subset X$ be an open stratum
of codimension $s$, let $\overline{S} \subset X$
be the closure, and let
$$\epsilon: \widetilde{S} \rightarrow X$$ be the normalization of
$\overline{S}$ equipped with a canonical map $\epsilon$
 to $X$. The normalization $\widetilde{S}$
is nonsingular and
separates the branches of the self-intersections of $\overline{S}$.
The map $\epsilon$ is an immersion locally on the source and therefore has a well-defined
normal bundle
$$\mathsf{N}_\epsilon = \epsilon^* T_X / T_{\widetilde{S}}\, $$
of rank $s$.

An open stratum $S\subset X$ of codimension $s$
is \'etale locally cut out by
$s$ branches of the full divisor $D$. These $s$ branches are partitioned
by monodromy orbits over $S$. Each monodromy orbit determines
a summand of $\mathsf{N}_\epsilon$. We  obtain a canonical splitting
of $\mathsf{N}_\epsilon$ corresponding to monodromy orbits
$$\mathsf{N}_\epsilon = \oplus_{\gamma\in \mathsf{Orb}(S)}
\mathsf{N}^\gamma_\epsilon\,,\  \ \ \ \mathsf{rank}(\mathsf{N}^\gamma_\epsilon)
= |\gamma|\, ,$$
where $\mathsf{Orb}(S)$ is the set of monodromy orbits of the
branches of $D$ cutting out $S$, and $|\gamma|$ is the number
of branches in the orbit $\gamma$.
For polynomials $P_\gamma$ in the Chern classes of $\mathsf{N}^\gamma_\epsilon$,
we define
\begin{equation}\label{dd7h}
[S, \{P_\gamma\}_{\gamma\in \mathsf{Orb}(S)}] =
\epsilon_*\left( \prod_{\gamma\in \mathsf{Orb}(S)} P_\gamma( \mathsf{N}^\gamma_\epsilon)\right)
\in \mathsf{CH}^\star(X)\, .
\end{equation}

We define {\em normally decorated classes} by the following more general construction. 
Let $G$ be the monodromy group of the $s$ branches of $D$ which cut out $S$. Over 
$\widetilde{S}$, there is a principal $G$-bundle
$$ \mu: \widetilde{P} \rightarrow \widetilde{S}$$
over which the $s$ branches determine $s$ line bundles 
\begin{equation}\label{llbb}
N_1,\ldots, N_s\, .
\end{equation}
The $G$-action on $\widetilde{P}$ permutes the line bundles \eqref{llbb}
via the original monodromy representation.
Let $P_G$ be any $G$-invariant polynomial in the Chern classes $c_1(N_i)$.
Since $P_G(c_1(N_1), \ldots,c_1(N_s))$ is $G$-invariant, 
$$P_G(c_1(N_1), \ldots,c_1(N_s)) \in \mathsf{CH}^\star(\widetilde{S})\, .$$
We define a {\em normally decorated strata class} by \Sam{maybe we should be more careful with the notation $\epsilon$ here, since we are now on the torsor}
$$[S, P_G] =\epsilon_*(P_G(c_1(N_1), \ldots,c_1(N_s)))
\in \mathsf{CH}^\star(X)\, .$$
Construction \eqref{dd7h} is a special case of a normally decorated strata class.

A fundamental result about the log Chow ring of $(X,D)$ is the following
inclusion.

\begin{theorem}\label{cc77}
  Let $(X,D)$ be a nonsingular variety with a normal crossings divisor.
  Let $S\subset X$ be an open stratum. Every normally decorated class
  associated to $S$ lies in $\mathsf{divlogCH}^\star(X,D)$.
\end{theorem}

{We give two proofs of Theorem \ref{cc77}: in Section \ref{lcr2} we give a very concrete iterated blow-up of $X$ and an explicit computation expressing the normally decorated class as a sum of products of divisors. On the other hand, in Corollary \ref{Cor:divlogCH} we give a more conceptual explanation based on the study of the Chow group of the Artin fan of the pair $(X,D)$. }

\subsection{Proof of Theorem \ref{cc77}} \label{lcr2}
Theorem \ref{cc77} is almost trivial if every irreducible component
$D_i$ of $D$ is nonsingular. The complexity of the argument occurs
only in case there are irreducible components with self-intersections.


\vspace{10pt}
\noindent{\em Proof.} 
 Let $S\subset X$ be an open stratum
of codimension $s$.
The first case to consider
is when $S$ is simple. Then, the closure
$$\overline{S} \subset X$$
in nonsingular and no normalization
is needed,
$$\epsilon: \overline{S} \rightarrow X\, .$$
Let $G$ be the monodromy of the $s$ branches of $D$
which cut out $S$.
We must prove
$$[S, P_G] =\epsilon_*(P_G(c_1(N_1), \ldots,c_1(N_s)))
\in \mathsf{divlogCH}^\star(X)\, $$
for every $G$-invariant polynomial $P_G$.

We argue by induction on the degree of $P_G$.
The base case is
when $P_G$ is
of degree 0. We
can take $P_G=1$, and we must prove
\begin{equation}
\label{vvpp}
[S, 1] =\epsilon_*[S]
\in \mathsf{divlogCH}^\star(X,D)\, . 
\end{equation}
Our argument
requires a
blow-up construction which we term 
an explosion.

The {\em explosion}
of $(X,D)$ along a simple stratum $S$,
\begin{equation}\label{crr4}
e:\mathsf{E}_S(X,D) \rightarrow X\,,
\end{equation}
is defined by a sequence of
blow-ups of $X$.
To describe the blow-ups locally{\footnote{Throughout the
proof of Theorem \ref{cc77}, the terms local, near, and
open refer to the Euclidean topology since we must 
separate branches.}}
near a point $p\in S$, let
$$B_1,\ldots,B_s$$
be the branches of $D$ cutting out $S$ near $p$.
\begin{enumerate}
\item[$\bullet$]
At the $0^{th}$ stage, we blow-up $S$, the intersection of all $s$ branches $B_1,\ldots,B_s$.
\end{enumerate}
Consider next the 
strict
transform of
the intersection
of  $s-1$ branches. For each choice of $s-1$ branches, the
strict transform of the intersection is
nonsingular of codimension $s-1$ over an open set of $p\in X$. Moreover,
the strict transforms of the intersections of different sets of $s-1$ branches
are disjoint over an open set of $p\in X$.
\begin{enumerate}
\item[$\bullet$]
At the $1^{st}$ stage, we blow-up all $s$ of these strict transforms of intersections of $s-1$
branches. 
\end{enumerate}
Then, the strict
transforms of
the intersections
of  $s-2$ branches among $B_1,\ldots, B_s$ are 
nonsingular of codimension $s-2$ and disjoint over an open set of $p\in X$.
\begin{enumerate}
\item[$\bullet$]
At the $2^{nd}$ stage,
we blow-up all $\binom{s}{2}$ of these strict transforms of intersections of $s-2$
branches. 
\end{enumerate}
We proceed in the above pattern until we have completed $s-1$ stages.
\begin{enumerate}
\item[$\bullet$]
At the $j^{th}$ stage, we blow-up all $\binom{s}{j}$
strict transforms
of intersections of $s-j$ branches. 
\end{enumerate}
The explosion \eqref{crr4} is
the result after stage
$j=s-1$.\footnote{At stage $j=s-1$, we are blowing-up divisors,
so no change occurs in the space. 
Still, to uniformize later notation, we include this $j=s-1$ stage and declare the divisorial center of this trivial blowup to be its exceptional divisor.
}
Since the above blow-ups are defined symmetrically with respect to the
branches $B_i$, the definition is well-defined globally on $X$.

Near $S$, all the prescribed blow-ups are of simple loci, but
non-simplicity may occur away from $S$. In order for the explosion to
be a sequence of simple blow-ups, some extra blow-ups may be
required
far from $S$. Since we will only be interested in the geometry near
$S$, the blow-ups related to non-simplicity away
from $S$ are not important for our argument (and are not
included in our notation).

A local study shows the following properties of the explosion
\begin{equation*}
e:\mathsf{E}_S(X,D) \rightarrow X\,,
\end{equation*}
near $S$:
\begin{enumerate}
    \item [(i)] The inverse image
    $ e^{-1}(S) \subset \mathsf{E}_S(X,D)$
is a nonsingular irreducible subvariety which we denote by $\mathsf{E}_S(S)$ and
call the {\em exceptional divisor} of the  explosion.
We denote the inclusion by
$$\iota: \mathsf{E}_S(S) \rightarrow \mathsf{E}_S(X,D)\, . $$
\item [(ii)] 
Let $\mathsf{N}_S$ be the rank $s$ normal bundle of $S$ in $X$.
The fibers of the projective normal bundle
\begin{equation}\label{xx44}
\mathsf{P}(\mathsf{N}_S) \rightarrow S
\end{equation}
have a canonical (unordered) set of $s$ coordinate hyperplanes
determined by the $s$ local branches of $D$ cutting out $S$.
In the fibers of \eqref{xx44},
these relative hyperplanes determine $s$ coordinate points,
$\binom{s}{2}$ coordinate lines, $\binom{s}{3}$ coordinate planes, and so on.
\item[(iii)]
The restriction of the explosion morphism to the exceptional divisor 
$$ e_S: \mathsf{E}_S(S) \rightarrow S$$
is obtained from $\mathsf{P}(\mathsf{N}_S) \rightarrow S$ by first blowing-up the coordinate points, and then
blowing-up the strict transforms of the coordinate lines, and so on.
For $$1\leq j \leq s-1\, ,$$
the $j^{th}$ stage of the construction of the
explosion restricts 
to the blow-up of the strict transform of the $(j-1)$-dimensional
coordinate linear spaces of the fibers of \eqref{xx44}.
\item[(iv)] 
On $\mathsf{E}_S(S)$, we have a distinguished set of divisors
$$E_0, E_1, \ldots, E_{s-1} \in \mathsf{CH}^1( \mathsf{E}_S(S))\, .$$
Here, $E_0$ is the pull-back to $\mathsf{E}_S(S)$ of
$$\mathcal{O}_{\mathsf{P}(\mathsf{N}_S)}(-1) \rightarrow \mathsf{P}(\mathsf{N}_S)$$
determined by the $0^{th}$ stage of the construction of the explosion.
Then, $E_j \in \mathsf{CH}^1( \mathsf{E}_S(S))$ is the pull-back to 
$\mathsf{E}_S(S)$ of the exceptional divisor obtained from the
blow-up of the strict transform of the $(j-1)$-dimensional coordinate
linear spaces in the fibers of \eqref{xx44}.
\item[(v)] Every class of the form
$$ [\mathsf{E}_S(S)]\cdot \mathsf{F}(E_0,\ldots,E_{s-1}) \in \mathsf{CH}^*(\mathsf{E}_S(X,D))\, $$
where $\mathsf{F}$ is a polynomial, lies in the divisor ring of log Chow,
$$ [\mathsf{E}_S(S)]\cdot \mathsf{F}(E_0,\ldots,E_{s-1}) \in \mathsf{divlogCH}^*(X,D)\, .$$
The claim follows from the geometric construction of the explosion.
To start, $\mathsf{E}_S(S)$ is a component of the associated normal crossings divisor of 
$\mathsf{E}_S(X,D)$. 
For each $0\leq j \leq s-1$,
$E_j$ comes from the pull-back of a divisor stratum of the blow-up at
the $j^{th}$ stage. 
\end{enumerate}

To the explosion geometry, we can apply
Fulton's excess intersection formula.  We start with the $0^{th}$ stage:
$$ e_0: X_0 \rightarrow X$$
is the blow-up along $S$, and
$$ e_0^*[S] = [\mathsf{P}(\mathsf{N}_S)] \cdot c_{s-1}\left(\frac{\mathsf{N}_S}{ \mathcal{O}_{\mathsf{P}(\mathsf{N}_S)}(-1)}\right)\, .$$
When we pull-back $e_0^*[S]$ all the way to $\mathsf{E}_S(X,D)$, we obtain{\footnote{
We have omitted the pull-backs in the notation inside the argument of $c_{s-1}$.}}
$$ e^*[S] = [\mathsf{E}_S(S)] 
\cdot c_{s-1}\left(\frac{\mathsf{N}_S}{ \mathcal{O}_{\mathsf{P}(\mathsf{N}_S)}(-1)}\right)\, .$$
By definition, we have
$$ c(\mathcal{O}_{\mathsf{P}(\mathsf{N}_S)}(-1)) = 1+ E_0\, .$$
By property (v) above for the explosion geometry, to
prove 
\begin{equation}\label{ckkd}
\epsilon_*[S]
\in \mathsf{divlogCH}^\star(X,D)\, ,
\end{equation}
we need only show
\begin{equation}\label{dd33}
c_{k}({\mathsf{N}_S}) = \mathsf{F}_k(E_0, \ldots, E_{s-1}) \in \mathsf{CH}^k(\mathsf{E}_S(S))
\end{equation}
for polynomials $\mathsf{F}_k$ for $1\leq k \leq s-1$.

The claim \eqref{dd33} is established directly by the following
basic formula of the explosion geometry. For $0\leq j\leq s-1$, let
$$\mathsf{L_j} = \sum_{i=0}^j E_i\, .$$
Let $\sigma_k$ be the $k^{th}$ elementary symmetric polynomial.
Then, {we claim that}
\begin{equation}\label{dd333}
c_{k}({\mathsf{N}_S}) = {\sigma}_k(\mathsf{L}_0, \ldots, \mathsf{L}_{s-1}) \in \mathsf{CH}^k(\mathsf{E}_S(S))\, .
\end{equation}
{
Once we prove \eqref{dd333}, this immediately shows \eqref{dd33} and thus, as explained above, establishes \eqref{ckkd}. We remind ourselves that \eqref{ckkd} represents the 
base case $P_G=1$ of our inductive proof that $[S, P_G] \in \mathsf{divlogCH}^*(X)$.}

Let $\mathsf{T}=(\com^*)^s$ and let
$t_i : T \to \com^*$ be the projection to the $i^{th}$ factor, which we interpret as the
weight of the standard representation of this $i^{th}$ factor.
To show formula \eqref{dd333}, we
consider the universal
$\mathsf{T}$-equivariant
model where $S\subset X$ is  
$$\mathbf{0}\in \com^s$$
and the logarithmic boundary $H\subset \com^s$
is the union of the $s$
coordinate hyperplanes.
Then, the $\mathsf{T}$-action on
$$e_{\mathbf{0}}:\mathsf{E}_{\mathbf{0}}(\com^s,H) \rightarrow {\mathbf{0}}$$
has $s!$ isolated $\mathsf{T}$-fixed points naturally indexed
by elements of the symmetric group $\Sigma_s$. 
\Joh{Is this correct to add: "Indeed, given $\gamma \in \Sigma_s$, there is a unique fixed point which in the $j$-th stage of the blow-up process lies over the strict transform of the intersection of branches $\{1, \ldots, s\} \setminus \{\gamma(1), \ldots, \gamma(j)\}$."}
The weights of the
divisors 
$$\mathsf{L}_0, \ldots, \mathsf{L}_{s-1}$$
with their canonical 
$\mathsf{T}$-equivariant
lifts at the $\mathsf{T}$-fixed point
$\gamma\in \Sigma_s$
are
$$t_{\gamma(1)},
t_{\gamma(2)}, t_{\gamma(3)}, \ldots ,
t_{\gamma(s)}\, $$
respectively.
Formula
\eqref{dd333} then follows
immediately for the
$\mathsf{T}$-equivariant
model. The general
case of \eqref{dd333}
is a formal consequence.

We now will establish
the induction step.
Let $S\subset X$
be a simple stratum of
codimension $s$
with monodromy 
group{\footnote{
The geometry
involved in the
proof of the base
case of the induction was
fully symmetric
with respect to the
branches, so the group
$G$ did not play a role.
}} $G$
of the branches of $D$
cutting out $S$.
We must prove
$$[S, P_G] =\epsilon_*(P_G(c_1(N_1), \ldots,c_1(N_s)))
\in \mathsf{divlogCH}^\star
(X,D)\, $$
for every $G$-invariant polynomial $P_G$. By induction, we  assume the truth of the
statement for polynomials
of lower degree.

Let $P_G$ be a $G$-equivariant polynomial in $c_1(N_1), \ldots, c_1(N_s)$
of degree $d>0$. 
We will prove a stronger property for the
induction argument:
$$\epsilon_*(P_G(c_1(N_1), \ldots,c_1(N_s)))
\in  \mathsf{divlogCH}^\star
(X,D)$$
{\em
can be expressed as
a linear combination of
terms of the form}
$$\widehat{D}_1\widehat{D}_2 \cdots \widehat{D}_d$$
{\em where the $\widehat{D}_i$
are components of the
logarithmic boundary of
an iterated blow-up of
the explosion 
$\mathsf{E}_S(X,D)$ and
$\widehat{D}_1$ lies over}
$$\mathsf{E}_S(S)\subset
\mathsf{E}_S(X,D)\, .$$ Our proof of the
base of the induction
 establishes the stronger
property.

We can assume $P_G$ is the summation{\footnote{The
stabilizer factor 
occurs to correct for
overcounting.}}
$M_G$ of the $G$-orbit of a degree $d$ monomial $M$,
$$M_G=
\frac{1}{|\mathsf{Stab}(M)|}
\sum_{g\in G} g(M)\, .
$$
We will study the geometry of the the exceptional divisor of the explosion
$$e_S: \mathsf{E}_S(S) \rightarrow S$$
locally over an analytic
open set 
$U_p \subset S$ of $p\in S$. 

Over small enough $U_p$, we can separate all the branches 
$B_1,\ldots, B_s$
of $D$ which cut out $S$,
and we can
 write
\begin{equation}\label{cctt5}
M = c_1(N_1)^{m_1}\cdots c_1(N_s)^{m_s} = B_1^{m_1} \cdots B_s^{m_s}\, .
\end{equation}
Over $U_p$, we can separate all the exceptional divisors of all
the blow-ups in the construction of
$$ \mathsf{E}_S(S) \rightarrow \mathsf{P}(\mathsf{N}_S)$$
explained in (iii) above. There are $2^{s}-2$ such exceptional divisor in bijective
correspondence to all the proper {nonzero} coordinate linear subspaces of the fiber
$\mathsf{N}_{S}|_p$
of $\mathsf{N}_S$ at $p$. We denote these $2^{s}-2$ exceptional divisors by
$E_\Lambda$ where $$\Lambda \subset \mathsf{N}_{S}|_p$$ is a proper coordinate linear space.
As before, we denote the pull-back of $\mathcal{O}_{\mathsf{P}(\mathsf{N}_S)}(-1)$
to $\mathsf{E}_S(S)$ by $E_0$.

Via the pull-back formula
for  $B_i$, we have 
\begin{equation}\label{kk339}
e^*(N_i) = E_0+ \sum_{\Lambda \subset H_i} E_\Lambda \  \in \mathsf{CH}^1( e^{-1}(U_p))\, ,
\end{equation}
where $H_i \subset \mathsf{N}_{S}|_p$  is the hyperplane associated to $B_i$.
We now substitute  formula \eqref{kk339} into \eqref{cctt5} to find
$$M\in \mathbb{Q}[E_0,
\{ E_\Lambda\}_\Lambda]\, .$$
\Joh{Maybe it would be better to write here something like $e^* M = \sum_{E \in \{E_0\} \cup \{ E_\Lambda\}_\Lambda} \mu_E \cdot  M^E \in \mathbb{Q}[E_0,
\{ E_\Lambda\}_\Lambda]$ above; otherwise I don't see a way to write the expression for $e^*[S, M_G]$ in a reasonable way below -- I think you need to have a linear combination with some coefficients $c_E$.}
Of course, 
$M$ has degree $d$ in
the divisors
$E_0$ and 
$\{ E_\Lambda\}_\Lambda$.

Let $M^E$ be a 
monomial of degree $d$ in the divisors 
\begin{equation}\label{dd99e}
E_0\ \ \text{and}\ \  
\{ E_\Lambda\}_\Lambda\, .
\end{equation}
The monodromy group $G$
acts{\footnote{
The $G$-action on $\{ E_\Lambda\}_\Lambda$
preserves the dimension of $\Lambda$. Moreover,
for a group element $\g\in G$, if $g(E_\Lambda) \neq E_\Lambda$, then
$$g(E_\Lambda) \cap E_\Lambda= \emptyset\, .$$
\label{ffff}
}} 
canonically on 
the set \eqref{dd99e} leaving $E_0$ fixed.
Let $$M^E_G=
\frac{1}{|\mathsf{Stab}(M^E)|}
\sum_{g\in G}
g(M^E)$$
be the summation over the
$G$-orbit of $M^E$.
Since $M^E_G$ is $G$-invariant, $M^E_G$
is a well-defined class
\begin{equation*}
M^E_G \in \mathsf{CH}^d(\mathsf{E}_S(S))\, .\end{equation*}

\Joh{Do we want to claim at this point that $e_S^* M$ can be written as a linear combination of the monomials $M_G^E$? I would say this is not entirely obvious since in the middle of our argument we used some equalities in Chow groups of "small enough $U_p$". Will there be correction terms outside of that open $U_p$?}

To prove the stronger
induction
step, we need only 
prove{\footnote{Recall, $\iota$ is the inclusion
$\iota: \mathsf{E}_S(S) \rightarrow \mathsf{E}_S(X,D)$.}}
\begin{equation}\label{kss9}
\iota_*\left(M^E_G
\cdot c_{s-1}\left(\frac{\mathsf{N}_S}{ \mathcal{O}_{\mathsf{P}(\mathsf{N}_S)}(-1)}\right)\right)
 \in \mathsf{divlogCH}^*(X,D)\,  \end{equation}
can be expressed as
a linear combination of
terms of the form
$$\widehat{D}_1\widehat{D}_2 \cdots \widehat{D}_d$$
where the $\widehat{D}_i$
are components of the
logarithmic boundary of
an iterated blow-up of
the explosion 
$\mathsf{E}_S(X,D)$ and
$\widehat{D}_1$ lies over
$\mathsf{E}_S(S)$.
To see why the claim for \eqref{kss9} is enough, we write
\begin{eqnarray*}
e^*[S,M_G]
& =& \sum_{M^E_G} e^*[S] \cdot M^E_G\\
&=& \sum_{M^E_G}
[\mathsf{E}_S(S)] 
\cdot c_{s-1}\left(\frac{\mathsf{N}_S}{ \mathcal{O}_{\mathsf{P}(\mathsf{N}_S)}(-1)}\right)
\cdot M^E_G
\\
& = & \sum_{M^E_G}
\iota_*\left(M^E_G
\cdot c_{s-1}\left(\frac{\mathsf{N}_S}{ \mathcal{O}_{\mathsf{P}(\mathsf{N}_S)}(-1)}\right)\right)
\, .
\end{eqnarray*}
\Joh{Referee: "Notation is unclear. What are sums over? How are the $M_G^E$ related to $M_G$?" I think currently these points are indeed a bit implicit.}
The first equality
is written with the
understanding that
$e^*[S]$ is supported
on $\mathsf{E}_S(S)$.
\Joh{Referee: "I do not understand the "with the understanding"-comment". I think we need this since some of the relations we use above only hold on the preimage of $S$ under $e$. But see also my comments above.}

To study $M^E_G$,
we take a geometric approach. If $M^E$
is just $E_0^d$, then 
\eqref{kss9} is already of the
claimed form by our
analysis in the base case.
Otherwise, $M^E$ has at
least one factor $E_\Lambda$.
Since $\{E_\Lambda\}_\Lambda$
is a set of simple normal crossings divisors
on $\mathsf{E}_S(S)$, we claim that we can write $M^E$ (if nonzero) 
as
$$M^E=E_{\Lambda_1}\cdots E_{\Lambda_t} \cdot \widetilde{M}^E\, ,$$
where $E_{\Lambda_1},\ldots,E_{\Lambda_t}$ are distinct divisors
with a nonempty transverse intersection
$$I_{U_p} = E_{\Lambda_1}\cap \ldots\cap E_{\Lambda_t}\  \ 
\text{over} \ \  U_p\, .$$
Moreover, we can assume
{\em every} divisor of the monomial
$\widetilde{M}^E$ contains $I_{U_p}$.
Indeed, we construct inductively for $i=1,2, \ldots$ a representation
$$
M^E = E_{\Lambda_1} \cdots E_{\Lambda_i} \cdot \widetilde{M}_i^E
$$
such that the $E_{\Lambda_j}$ are distinct and have nonzero, transverse intersection. For $i=1$ this is just our assumption that $M^E$ has some factor $E_{\Lambda}=:E_{\Lambda_1}$. On the other hand, given the representation above for some $i$, if all factors $E_{\Lambda'}$ of $\widetilde{M}_i^E$ contain $E_{\Lambda_1} \cap \ldots \cap  E_{\Lambda_i}$, we are done, setting $t=i$. If there is an $E_{\Lambda'}$ not satisfying this, we set $E_{\Lambda_{i+1}}=E_{\Lambda'}$. If the intersection $E_{\Lambda_1} \cap \ldots \cap  E_{\Lambda_{i+1}}$ was empty, then $M^E=0$, giving a contradiction. Thus the intersection is nonempty, and transverse by the fact that the $E_{\Lambda}$ are a normal crossings divisor. We continue inductively and this construction concludes after at most $d$ steps.

When the monodromy invariant $M^E_G$ is
considered, we obtain a nonsingular  subvariety
of  $\mathsf{E}_S(S)$
of codimension $t$,
\Joh{Referee: "How do we obtain it? What do we obtain it from?"}
$$ V \subset \mathsf{E}_S(S)$$
which is a simple stratum of $\mathsf{E}_S(X,D)$,
$$\epsilon^V: V \rightarrow \mathsf{E}_S(X,D)\, .$$
Over $U_p$, the subvariety $V$ restricts
to the union{\footnote{The distinct $G$-translates of $I_{U_p}$ are disjoint, see Footnote \ref{ffff}.}} of
the distinct
$G$-translates of $I_{U_p}$.
The crucial geometric observation is  
$$\iota_*(M^E_G) = \epsilon^V_*(\widetilde{P}) \in \mathsf{CH}^\star(\mathsf{E}_S(X,D)\, ,$$
where $\widetilde{P}$ is defined by $\widetilde{M}^E$
and is of degree at most $d-1$.

We can apply the strong induction property: the class
$$\epsilon^V_*(\widetilde{P})
\in  \mathsf{divlogCH}^\star
(X,D)$$
can be expressed as
a linear combination of
terms of the form
$$\widehat{D}_1\widehat{D}_2 \cdots \widehat{D}_d$$
where the $\widehat{D}_i$
are components of the
logarithmic boundary of
an iterated blow-up of
the explosion of $V$ in $\mathsf{E}_S(X,D)$
and
$\widehat{D}_1$ lies over
$$\mathsf{E}_V(V)\subset
\mathsf{E}_S(X,D)\, .$$
Then, the claim
\begin{equation}\label{nn55p}
\iota_*\left(M^E_G
\cdot c_{s-1}\left(\frac{\mathsf{N}_S}{ \mathcal{O}_{\mathsf{P}(\mathsf{N}_S)}(-1)}\right)\right)
 \in \mathsf{divlogCH}^*(X,D)\,  \end{equation}
holds by the analysis of
$$c_{s-1}\left(\frac{\mathsf{N}_S}{ \mathcal{O}_{\mathsf{P}(\mathsf{N}_S)}(-1)}\right)$$
on $\mathsf{E}_S(S)$ in the base case of the induction.
Since each monomial
$$\widehat{D}_1\widehat{D}_2 \cdots \widehat{D}_d$$
of $\epsilon^V_*(\widetilde{P})$ lies over 
$\mathsf{E}_V(V)$ which,
in turn, lies over $\mathsf{E}_{S}(S)$,
the analysis of the base case yields the
desired result \eqref{nn55p}.

The induction argument is complete, so we have
proven Theorem \ref{cc77} in case $S$ is a simple
stratum of $(X,D)$.  The general case follows
by repeated application of the result for a simple
stratum. 

Let $S\subset X$ be a stratum with
a singular closure
$$\overline{S} \subset X\, .$$
The first step
is to blow-up simple strata in $\overline{S}$,
$$\widehat{X} \rightarrow X\, ,$$
until the strict transform of $\overline{S}$, 
$$\widehat{S} \subset \widehat{X}\, ,$$
is nonsingular.
Since $S$ is simple stratum of the blow-up $\widehat{X}$,
we can apply Theorem \ref{cc77} to ${S} \subset \widehat{X}$.

Via the blow-down map, we have
$$ \widehat{S} \rightarrow \overline{S}\, . $$
There are two discrepancies
to handle before deducing Theorem \ref{cc77} for
normally decorated classes associated to $S\subset X$
from the result for normally decorated classes
associated to $S\subset \widehat{X}$:
\begin{enumerate}
    \item [(i)]
    The fundamental class $[\widehat{S}]\in \mathsf{CH}^\star(\widehat{X})$ is not the
    pull-back of $[\overline{S}]\in \mathsf{CH}^\star({X})$.
    \item[(ii)] The normal directions of 
    $\widehat{S}\subset \widehat{X}$ differ
    from the pull-backs of the normal directions
    of $\overline{S}\subset {X}$.
    \end{enumerate}
However, both discrepancies are corrected by
applying the simple stratum result to the
lower dimensional strata occurring in 
$\widehat{S} \setminus S$.
\qed
\vspace{10pt}

\subsection{Explosion geometry and barycentric subdivision} \label{expbar}
The explosion operation $E(X,D)$ along a simple stratum $S\subset X$, which appeared in the proof \ref{lcr2}, is an essentially  combinatorial operation that has a natural interpretation in terms of the geometry of the cone complex $C(X,D)$.

 Consider first a cone $\sigma$ of dimension $n$ in a lattice $N$, and let $A_\sigma$ be the associated toric variety. Let $\mathcal{A}_\sigma$ be the
 associated Artin fan, which is simply the stack quotient of $A_\sigma$ by the corresponding dense torus $T_\sigma$. The logarithmic stratification of $A_\sigma$ is precisely the stratification defined by the orbits of $T_\sigma$, and there is a bijective dimension reversing correspondence between faces of $\sigma$ and strata. We write $\sigma(k)$ for the $k$-dimensional faces of $\sigma$ and thus the codimension $k$ strata of $A_\sigma$. 

For each face $\tau$ of $\sigma$, the barycenter $b_\tau$ of $\tau$ is the sum $$b_\tau=
\sum_{v_i \in \tau \cap \sigma(1)} v_i$$ of the primitive vectors along the extremal rays of $\tau$. For any flag $$\tau_0 \subset \tau_1 \subset \cdots \subset \tau_k$$ of faces of $\sigma$, the barycenters $b_{\tau_0}, \cdots, b_{\tau_k}$ span a cone. The set of all such cones, for all flags in $\sigma$, forms a subdivision of $\sigma$, which we call the \emph{barycentric subdivision} $\widetilde{\sigma}$ of $\sigma$. 

Alternatively, we can build the barycentric subdivision inductively: at step $1$, we start with the star subdivision over the barycenter of faces in $\sigma(n)$ (where $\sigma$ has dimension $n$), then take the star subdivision over faces in $\sigma(n-1)$, and so on, terminating after $n-1$ steps with $\sigma(2)$, after which the operation no longer has effect. We thus produce a sequence of $n-1$ subdivisions $$\widetilde{\sigma} = \sigma_{n-1} \to \sigma_{n-1} \cdots \to \sigma_1 \to \sigma_0=\sigma $$  When $\sigma = \mathbb{R}_{\ge 0}^n$, which is our main case of interest, the barycentric subdivision has $n!$ maximal cones. 

The barycentric subdivision of $\sigma$ produces a log modification $$\widetilde{A}_\sigma \to A_\sigma\, ,$$ which is in fact a log blow-up. More precisely, we have constructed the subdivision $\widetilde{A}_\sigma \to A_\sigma$ as a sequence 
$$
\widetilde{A}_\sigma = A_{n-1} \to A_{n-2} \to \cdots A_1 \to A_0 = A_\sigma
$$

\noindent and the map $A_k \to A_{k-1}$ is determined by the subdivision $\sigma_k \to \sigma_{k-1}$, which is the subdivision corresponding to the domains of linearity of a piecewise linear function -- see \cite{KKMSD} for the construction. In the case of interest, $$\sigma = \mathbb{R}_{\ge 0}^n\, ,$$ the map $A_1 \to A_0$ is the blowup of $\mathbb{A}^n$ at the origin, $A_2 \to A_1$ is the blowup along the strict transforms of the coordinate lines, and in general $A_k \to A_{k-1}$ is the blowup along the strict transforms of the dimension $k-1$ hyperplanes of $\mathbb{A}^n$ in $A_{k-1}$. Thus, the barycentric subdivision of $\mathbb{A}^n$ is precisely the explosion of $\mathbb{A}^n$ along the origin.   

The barycentric subdivision construction is clearly equivariant and therefore descends to the Artin fan $\mathcal{A}_\sigma$ of $A_\sigma$. Furthermore, the
subdivision is the same on  isomorphic faces of $\sigma$ and invariant with respect to  automorphisms of $\sigma$. Consequently, given any cone complex $C$, the barycentric subdivisions of individual cones glue to a global subdivision of $C$, and that is true even if faces of $C$ are identified or if there is monodromy in $C$. Thus, for a normal crossings pair $(X,D)$, we can define the barycentric subdivision $\widetilde{C}(X,D)$ of the cone complex $C(X,D)$, and equivalently, a log blow-up $$\widetilde{\mathcal{A}}_X\to
\mathcal{A}_X$$
of the Artin fan.
We also obtain \emph{globally} a log blow-up $$(\widetilde{X},\widetilde{D})= X \times_{\mathcal{A}_X} \widetilde{\mathcal{A}}_X \to (X,D)$$  with Artin fan $\mathcal{A}_{\widetilde{X}} = \widetilde{\mathcal{A}}_X$. 

The explosion of Section \ref{lcr2} can only be defined locally around a simple stratum $S$. 
A quasi-projective stratum $S$ (not necessarily simple) of a normal crossings pair $(X,D)$ corresponds to a cone $\sigma$ of $C(X,D)$.  More precisely, the quasi-projective stratum $S$ corresponds to the interior of $\sigma$, and the whole of $\sigma$ corresponds to a canonical open set $U$ in $X$ that contains $S$ as its minimal stratum: the open set $U$ consists of all quasi-projective strata whose closure contains $S$. The explosion
$\mathsf{E}_S(U, D|_U))$ is well-defined.

The cone $\sigma$ has a cover by $\mathbb{R}^n_{\ge 0}$, and, more precisely, by a quotient of $\mathbb{R}^n_{\ge 0}$ obtained by potentially identifying faces and taking a quotient by a group $G$. The group $G$ is precisely the monodromy group of the divisors $D$ that cut out $S$ considered in Section \ref{dslc}, and the interior $\sigma^\circ$ of $\sigma$ is in fact the stack quotient $[\mathbb{R}_{>0}^n/G]$. Similarly, the Artin fan $\mathcal{U}$ of $U$ has an analogous \'etale cover by the groupoid quotient of $[\mathbb{A}^n/\mathbb{G}_m^n \rtimes G]$, with $S$ corresponding to the minimal stratum $$B(\mathbb{G}_m^n \rtimes G) \subset \mathcal{U}\, .$$
The cover is not representable, but is representable over $S$.  
From the discussion of the barycentric subdivision of $\mathbb{A}^n$, we see that $\mathsf{E}_S(U,D_U)$  is precisely the barycentric subdivision $\widetilde{X} \to X$ restricted to $U$. We may thus view the barycentric subdivision as globalizing the explosion geometry. 

If the stratum $S$ is simple, the explosion  of Section \ref{lcr2}  is
defined over a neighborhood of $\overline{S}$. However, the extension no longer coincides with the barycentric subdivision. The barycentric subdivision performs additional blowups, first blowing up all minimal strata in the closure of $S$ (and also strata around $\overline{S}$ whose closure does not necessarily meet $S$). 
\Joh{Maybe add: "Thus the explosion corresponds to the barycentric subdivision of the cone corresponding to $S$ and its faces (in decreasing order of combinatorial dimension).}
\Joh{Referee: "It looks like the explosion is the barycentric subdivision of the cone corresponding to $S$ and its faces (in decreasing order of combinatorial dimension)" Is that right? If so, we could add it?} \Sam{This is correct, and we say it in the sentence "From the discussion of the barycentric...". Maybe we can just make the sentence clearer?}
\Joh{I made a proposal based on the referee comment, feel free to modify.}
\\

We illustrate the concepts discussed above through an example. Let $(X,D)$ be a log scheme whose cone complex is the cone over an equilateral triangle, with all edges identified and with monodromy $\mathbb{Z}/3\mathbb{Z}$. For example, we can construct $(X,D)$  by taking $$X \to B$$ to be a family with fiber $\mathbb{A}^3$ over a nonsingular base $B$ satisfying $\pi_1(B) = \mathbb{Z}$, so that the generator of $\pi_1(B)$ cyclically
permutes the coordinate hyperplanes of
$\mathbb{A}^3$. The divisor $D\subset X$
is then the union of these coordinate
hyperplanes over $B$.

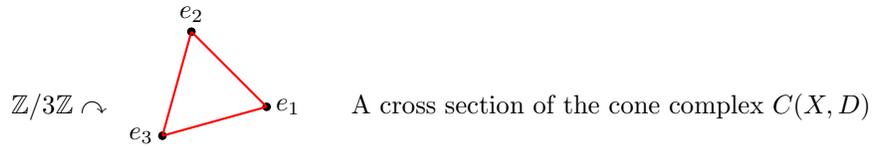
$$
  \begin{tikzpicture}
    \node[left] at (-1,0,0){$\mathbb{Z}/3\mathbb{Z} \curvearrowright$};
    \node[draw,circle,inner sep=1pt,fill] at (1, 0,0){};
    \node[right] at (1,0,0){$e_1$};
    \node[draw,circle,inner sep=1pt,fill] at (0, 1,0){};
    \node[above] at (0,1,0){$e_2$};
    \node[draw,circle,inner sep=1pt,fill] at (0, 0,1){};
    \node[left] at (0,0,1){$e_3$};
     \node[right] at (2,0,0){A cross section of the cone complex $C(X,D)$};
    
\draw[-, red, thick] (1,0,0) -- (0,1,0);
\draw[-, red, thick](0,1,0) -- (0,0,1);
\draw[-, red, thick](0,0,1) -- (1,0,0);

\end{tikzpicture}
$$

The log scheme $(X,D)$  has four strata: the open set $X-D$, corresponding to the empty face of the triangle (or, equivalently, the vertex of the cone over the triangle), the interior of the divisor $D$ corresponding to the vertex $$e_1=e_2=e_3\, ,$$ the locus which is \'etale locally the intersection of exactly two irreducible components of $D$ corresponding to edge
$$\overline{e_1e_2} = \overline{e_1e_3}= \overline{e_2e_3}\, ,$$ and the triple point singularity of $D$ corresponding to the whole triangle. 
We name the strata $Q,R,S,T$ respectively. While $T$ is simple, $S$ is not, since $$\overline{S}=S \cup T$$ is not normal. The strata are taken bijectively to points of the Artin fan via the map $$\alpha: X \to \mathcal{A}_X\, $$ 
We depict the Artin fan as four points, each isomorphic to $B\mathbb{G}_m^k \rtimes G$ as indicated, with points drawn increasingly bigger to describe the topology (the closure contains all smaller points). 

$$
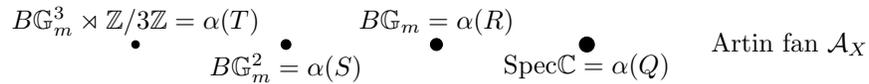
\begin{tikzpicture}
\node[draw,circle,inner sep = 1pt, fill] at (-1,0){};
\node[draw, circle, inner sep = 1.3pt, fill] at (1,0){};
\node[draw,circle,inner sep = 1.6pt, fill] at (3,0){};
\node[draw,circle,inner sep = 2pt, fill] at (5,0){};
\node[above] at (-1,0){$B\mathbb{G}_m^3 \rtimes \mathbb{Z}/3\mathbb{Z} = \alpha(T)$};
\node[below] at (1,0){$B\mathbb{G}_m^2 = \alpha(S)$}; 
\node[above] at (3,0){$B\mathbb{G}_m = \alpha(R)$}; 
\node[below] at (5,0){$\textup{Spec} \mathbb{C} = \alpha(Q)$};
\node[right] at (6,0){\, \, \, Artin fan $\mathcal{A}_X$};
\end{tikzpicture}
$$

Consider the explosion of the
quasi-projective stratum $S$
depicted by the open
line segment $\overline{e_1e_2}$. The open set $U$ over which the explosion is defined is $Q \cup R \cup S$.  The explosion of $S$ is the barycentric subdivision of 
$\overline{e_1e_2}$:  
$$
   
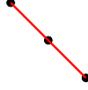
\begin{tikzpicture}
    \node[draw,circle,inner sep=1pt,fill] at (1, 0,0){};
    \node[right] at (1,0,0){};
    \node[draw,circle,inner sep=1pt,fill] at (0, 1,0){};
    \node[above] at (0,1,0){};
    \node[draw,circle, inner sep=1pt,fill] at (1/2, 1/2,0){};
    \draw[-,red, thick] (1,0,0) -- (0,1,0);

\end{tikzpicture}
$$
However, the above explosion does not extend away from $U$. The blowup of $\overline{S}$, over an \'etale cover of $X$ is depicted as 
$$
    
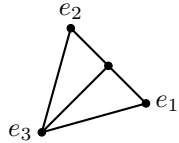
\begin{tikzpicture}
    \node[draw,circle,inner sep=1pt,fill] at (1, 0,0){};
    \node[right] at (1,0,0){$e_1$};
    \node[draw,circle,inner sep=1pt,fill] at (0, 1,0){};
    \node[above] at (0,1,0){$e_2$};
    \node[draw,circle,inner sep=1pt,fill] at (0, 0,1){};
    \node[left] at (0,0,1){$e_3$};
    \node[draw,circle,inner sep=1pt,fill] at (1/2, 1/2,0){};

\draw[-, thick] (1,0,0) -- (0,1,0);
\draw[-,thick] (1,0,0) -- (0,0,1); 
\draw[-, thick](0,1,0) -- (0,0,1); 
\draw[-,thick](0,0,1) --(1/2,1/2,0);

\end{tikzpicture}
$$
But the blow-up does not descend to $X$ as it does not respect the face identifications/automorphisms of $C(X,D)$. The barycentric subdivision is depicted as
$$
    
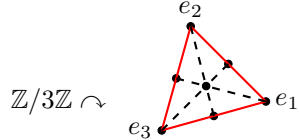
\begin{tikzpicture}
    \node[left] at (-1,0,0){$\mathbb{Z}/3\mathbb{Z} \curvearrowright$};
    \node[draw,circle,inner sep=1pt,fill] at (1, 0,0){};
    \node[right] at (1,0,0){$e_1$};
    \node[draw,circle,inner sep=1pt,fill] at (0, 1,0){};
    \node[above] at (0,1,0){$e_2$};
    \node[draw,circle,inner sep=1pt,fill] at (0, 0,1){};
    \node[left] at (0,0,1){$e_3$};
    \node[draw,circle,inner sep=1pt,fill] at (1/2, 1/2,0){};
    \node[draw,circle,inner sep=1pt,fill] at (1/2, 0,1/2){};
    \node[draw,circle,inner sep=1pt,fill] at (0, 1/2,1/2){};
    \node[draw,circle,inner sep=1pt,fill] at (1/3, 1/3,1/3){};

\draw[-, red, thick] (1,0,0) -- (0,1,0);
\draw[-, red, thick](0,1,0) -- (0,0,1);
\draw[-, red, thick](0,0,1) -- (1,0,0);
\draw[-, dashed, thick](1,0,0) -- (1/3,1/3,1/3);
\draw[-, dashed, thick](0,1,0) -- (1/3,1/3,1/3);
\draw[-, dashed, thick](0,0,1) -- (1/3,1/3,1/3);
\draw[-, dashed, thick](1/2,1/2,0) -- (1/3,1/3,1/3);
\draw[-, dashed, thick](1/2,0,1/2) -- (1/3,1/3,1/3);
\draw[-, dashed, thick](0,1/2,1/2) -- (1/3,1/3,1/3);

        \end{tikzpicture}
$$
The corresponding log blow-up restricts to the explosion over $U$. Over $X$,
the log blow-up is not the blow-up of $\overline{S}$, but the explosion of $T$.

\subsection{Tautological classes}\label{tautc}
Let $(X,D)$ be a nonsingular variety with a normal crossings divisor. We define the 
{\em logarithmic tautological ring}
$$\mathsf{R}^\star(X,D) \subset \mathsf{CH}^\star(X)$$
to be the $\mathbb{Q}$-linear subspace spanned by all 
normally decorated strata classes (which is
easily seen to be closed under the intersection product).
Theorem \ref{cc77} can then be written as
$$\mathsf{R}^\star(X,D) \subset 
 \mathsf{divlogCH}^\star(X,D)\, .$$

The logarithmic tautological ring
of $(X,D)$ depends strongly
on the divisor $D$. For example,
if $X$ is irreducible and $D=\emptyset$, then there is
only one stratum and
$$\mathsf{R}^\star(X,\emptyset) = \mathbb{Q}\, .$$
For the moduli space of curves, 
the inclusion
$$\mathsf{R}^\star(\overline{\mathcal{M}}_g, \Delta_0) \subset
\mathsf{R}^\star(\overline{\mathcal{M}}_g, \partial \overline{\mathcal{M}}_g)\, ,$$
is 
proper for $g\geq 2$. Furthermore,
the inclusion
$$\mathsf{R}^\star(\overline{\mathcal{M}}_g, \partial \overline{\mathcal{M}}_g) \subset
\mathsf{R}^\star(\overline{\mathcal{M}}_g)$$
in the standard tautological 
ring{\footnote{$\mathsf{R}^\star(\overline{\mathcal{M}}_g)$ is definitely
not equal to $\mathsf{R}^\star(\overline{\mathcal{M}}_g,\emptyset)$!
}}
is proper for $g\geq 3$
since
$\mathsf{R}^\star(\overline{\mathcal{M}}_g)$
contains $\kappa$ and $\psi$ classes
which do not appear in the
logarithmic constructions.

Let $(X,D)$ be a nonsingular variety with a normal crossings divisor.
Let $$\pi: \widetilde{X} \rightarrow X$$ be a simple blow-up
of $(X,D)$. Let $\widetilde{D}\subset \widetilde{X}$ be the
associated normal crossings divisor. We will prove the following two basic properties of 
logarithmic tautological rings.

\begin{theorem}\label{plbk}
The pull-back
$$\pi^* :
\mathsf{R}^\star({X},{D}) \rightarrow
\mathsf{CH}^\star(\widetilde{X})$$
has image in $\mathsf{R}^\star(\widetilde{X},
\widetilde{D})$.
\end{theorem}

\begin{theorem}\label{pfwd}
The push-forward
$$\pi_* :\mathsf{R}^\star(\widetilde{X},\widetilde{D}) \rightarrow
\mathsf{CH}^\star(X)$$
has image in $\mathsf{R}^\star(X,D)$.
\end{theorem}

By Theorems \ref{plbk} and \ref{pfwd}, we can simply write
$$\pi^* :\mathsf{R}^\star({X},{D}) \rightarrow
\mathsf{R}^\star(\widetilde{X},\widetilde{D})\,, \ \ \ 
\pi_* :\mathsf{R}^\star(\widetilde{X},\widetilde{D}) \rightarrow
\mathsf{R}^\star(X,D)\, .$$
 Theorems \ref{plbk} and \ref{pfwd}
 will proven in Section \ref{craf2} 
 via the geometry of the
 Artin fan. As a consequence,
 we will present
a more conceptual (but less constructive) proof of Theorem \ref{cc77}.

\subsection{The Chow ring of the
Artin fan} \label{craf}
Let $(X,D)$ be a nonsingular variety with a normal crossings divisor.
We relate here the normally decorated strata classes of $(X,D)$ to Chow classes on the Artin fan $\mathcal{A}_X$ of $(X,D)$.
Here, since $\mathcal{A}_X$ is a smooth, finite type algebraic stack stratified by quotient stacks, it has well-defined Chow groups $\mathsf{CH}^\star(\mathcal{A}_X)$ with an intersection product as defined in \cite{kreschcycleartin}.
{
Note that for our proof below it will not be necessary to recall the precise definition from \cite{kreschcycleartin}, since we only use some properties and examples of these Chow groups (like the existence of an excision sequence) that we recall when needed. 
Also, we stress again that all Chow groups below are with $\mathbb{Q}$-coefficients.
}

As we explain in
Section \ref{afan}, there is a smooth morphism
to the Artin fan,
$$\alpha: X \rightarrow \mathcal{A}_X\, .$$ 

\begin{theorem}\label{gg77}
There is a canonical isomorphism $$\mathsf{CH}^\star(\mathcal{A}_X) \cong \textup{PP}^\star(C(X,D))$$ between the Chow ring of $\mathcal{A}_X$ and
the algebra
of piecewise polynomial functions on the cone complex $C(X,D)$. 
\end{theorem}

\noindent{\em Proof.} By construction, the Artin fan $\mathcal{A}_X$ has a presentation as a colimit 
$$\mathcal{A}_X =
\varinjlim_{x \in \mathcal{S}} \mathcal{A}_{x}\, ,$$
where $\mathcal{S}$ is a finite diagram, each map $\mathcal{A}_x$ is a stack of the form $[\mathbb{A}^n/\mathbb{G}_m^n]$, and all maps in the diagram are \'etale. 
First, we note that for the individual stacks $\mathcal{A}_x = [\mathbb{A}^n/\mathbb{G}_m^n]$ we have
\begin{equation} \label{eqn:polyonsinglecone}
\mathsf{CH}^\star([\mathbb{A}^n/\mathbb{G}_m^n]) \cong \mathsf{CH}^\star([\mathsf{Spec}(\mathbb{C})/\mathbb{G}_m^n]) \cong  \mathbb{Q}[x_1, \ldots, x_n]\,.
\end{equation}
The first equality is because $$[\mathbb{A}^n/\mathbb{G}_m^n] \to [\mathsf{Spec}(\mathbb{C})/\mathbb{G}_m^n]$$ is a vector bundle and induces an isomorphism of Chow groups by \cite[Theorem 2.1.12 (vi)]{kreschcycleartin}. The second equality is
because the equivariant Chow ring of a product of tori is a polynomial algebra \cite[Section 3.2]{edidingraham}, which can be  identified with polynomials on the cone $\sigma_{X,x}$ associated to $\mathcal{A}_x$ (appearing in the colimit presentation of $C(X,D)$). 


For the entire Artin fan $\mathcal{A}_X$, we claim 
\begin{equation}\label{ddll}
\mathsf{CH}^\star{\mathcal{A}_X} = \varprojlim_{x \in \mathcal{S}} \mathsf{CH}^\star{\mathcal{A}_{x}}\, .
\end{equation}
If we can show equality \eqref{ddll}, then Theorem \ref{gg77} will follow since the result holds for each term on the right hand side by \eqref{eqn:polyonsinglecone}. Piecewise polynomial functions on $C(X,D)$ are defined by the corresponding limit presentation.

All the stacks appearing in \eqref{ddll} are very special: they are nonsingular and have a stratification with strata isomorphic to $$B(\mathbb{G}_m^n \rtimes G)$$ with $G$ a finite group. For the argument below, it will be more convenient to index Chow groups by the dimension of the cycles (instead of the codimension) and prove\footnote{A similar formula and computation for the Chow groups of the stack of expanded pairs appears in \cite{oesinghaus}.}
\begin{equation}\label{ccll}
\mathsf{CH}_{\star}(\mathcal{A}_X) = \varprojlim_{x \in \mathcal{S}} \mathsf{CH}_{\star}(\mathcal{A}_x)\, .
\end{equation}

Let $\mathcal{C}$ 
denote the full 2-subcategory of the 2-category of algebraic stacks with Ob($\mathcal{C}$) given by algebraic stacks $\mathcal{A}$ with a stratification by stacks of the form $B (\mathbb{G}_m^n \rtimes G)$, with $G$ a finite group.
Similarly, let $\mathcal{C}^\circ$
be the full 2-subcategory of $\mathcal{C}$ with objects given
by stacks of the form $B \mathbb{G}_m^n$. We start with a stack\footnote{The case of interest is the Artin fan $\mathcal{A}_X$ of $X$, but in the argument we allow $\mathcal{A}_X$ to be arbitrary in $\mathcal{C}$ in order to run the induction.}
$\mathcal{A}_X \in \mathcal{C}$  with a colimit presentation $$\mathcal{A}_X = \varinjlim_{x \in \mathcal{S}} \mathcal{A}_x = \mathcal{A}_X$$ where  $\mathcal{A}_x \in \mathcal{C}^\circ$ and all maps in the diagram are \'etale. We will prove \eqref{ccll} by induction on the number of strata of $\mathcal{A}_X$. 


Assume first that there is a unique stratum, $$\mathcal{A}_X = B(\mathbb{G}_m^n \rtimes G)\, ,$$ and all maps in the diagram $\mathcal{S}$ are isomorphisms. Then the groupoid $$\varinjlim_{x \in \mathcal{S}} \mathcal{A}_x$$ is equivalent to the quotient $B\mathbb{G}_m^n/G$, and the statement is equivalent to $$\mathsf{CH}_*(B(\mathbb{G}_m^n \rtimes G)) = \mathsf{CH}_*(B\mathbb{G}_m^n)^G\,,$$
which is true (see \cite[Lemma 2.20]{BaeSchmitt}). In general, we pick an open stratum $U \in \mathcal{A}_X$ with preimage $U_x \in \mathcal{A}_x$. Then, by \cite[Proposition 4.2.1]{kreschcycleartin} we have an exact sequence
\[
\begin{tikzcd}
\mathsf{CH}(U,1) \ar[r] & \mathsf{CH}(Z) \ar[r] & \mathsf{CH}(\mathcal{A}_X) \ar[r] & \mathsf{CH}(U) \ar[r] & 0 
\end{tikzcd}
\]  
with $Z = \mathcal{A}_X - U$. Since $U$ is of the form $U = B(\mathbb{G}_m^n \rtimes G)$, we can use \cite[Proposition 2.14, Remark 2.21]{BaeSchmitt} to see that
\[
\mathsf{CH}(U,1) = \mathsf{CH}(U) \otimes_{\mathbb Q} \mathsf{CH}(\mathrm{Spec}(\mathbb{C}), 1)
\]
Then by \cite[Remark 2.18]{BaeSchmitt}, the connecting homomorphism $\mathsf{CH}(U,1) \to \mathsf{CH}(Z)$ vanishes. So we obtain an exact sequence
\[
\begin{tikzcd}
0 \ar[r] & \mathsf{CH}(Z) \ar[r] & \mathsf{CH}(\mathcal{A}_X) \ar[r] & \mathsf{CH}(U) \ar[r] & 0\, , \\
\end{tikzcd}
\]  
and the same sequence holds with $\mathcal{A}_X$ replaced by $\mathcal{A}_x$, $U$ by $U_x$, and $Z$ by $Z_x = \mathcal{A}_x - U_x$. As projective limits are left exact, we obtain a diagram

\[
\begin{tikzcd}
0 \ar[r] & \mathsf{CH}(Z) \ar[r] \ar[d] & \mathsf{CH}(\mathcal{A}_X) \ar[r] \ar[d] & \mathsf{CH}(U) \ar[r] \ar[d] & 0 \\
0 \ar[r] & \varprojlim_{x \in \mathcal{S}}\mathsf{CH}(Z_x) \ar[r] & \varprojlim_{x \in \mathcal{S}}\mathsf{CH}(\mathcal{A}_x) \ar[r] & \varprojlim_{x \in \mathcal{S}}\mathsf{CH}(U_x)  &  \\
\end{tikzcd}
\]  

By induction, the left and right vertical arrows are isomorphisms. But the bottom row is exact as well: the composed map $$\mathsf{CH}(\mathcal{A}_X) \to \mathsf{CH}({U}) \cong \varprojlim_{x \in \mathcal{S}}\mathsf{CH}(U_x)$$ is surjective and factors through $\varprojlim_{x \in \mathcal{S}}\mathsf{CH}(\mathcal{A}_x)$. Thus the map $$\mathsf{CH}(\mathcal{A}_X) \to \varprojlim_{x \in \mathcal{S}}\mathsf{CH}(\mathcal{A}_x)$$ is an isomorphism as well. \qed

\vspace{8pt}
Theorem \ref{gg77} has clear
precursors in the toric
context by Payne \cite{Payne} and
Brion \cite{BrionPP}. 
In the logarithmic context, we were
directly motivated by
ideas of Ranganathan. A development of the theory for general log schemes will appear in \cite{MoRa}.

\begin{theorem}
\label{tautological=pp2}
The logarithmic tautological ring $$\mathsf{R}^\star(X,D) \subset \mathsf{CH}^\star(X)$$ coincides with the image $\alpha^*\mathsf{CH}^\star(\mathcal{A}_X) \subset \mathsf{CH}^\star(X)$.
\end{theorem}

\vspace{8pt}

\noindent {\em Proof.} \Sam{Newest version. I was trying to fix the old version and came up with this argument. It seems much simpler. Can you let me know if you believe it?}
Fix a stratum $S\subset X$ with closure $\overline{S}\subset X$, and  normalization 
$$\epsilon: \widetilde{S} \rightarrow \overline{S} \subset X\, .$$
Consider the cone complex $C(X,D)$  and the Artin fan $\mathcal{A}_X$ of $(X,D)$ with 
$$\alpha:X \to \mathcal{A}_X\, .$$ 

Let $\widetilde{P}$ be the total space of the principal $G$-bundle 
over the normalization $\widetilde{S}$ defined by the branches of $D$ in
Section \ref{deflc},
$$\mu: \widetilde{P} \to \widetilde{S}\,, \ \ \  \mu_X = \epsilon \circ \mu : \widetilde{P} \to X\, .$$
We observe that all the relevant geometry is pulled back from the Artin fan $\mathcal{A}_X$: 
the stratum $S$ corresponds to the stratum $$\alpha(S) = \mathcal{S} \subset \mathcal{A}_X\, $$
with closure $\overline{\mathcal{S}} = \alpha(\overline{S})$. Let $\widetilde{\mathcal{S}}$ be the normalization of $\overline{\mathcal{S}}$, and let 
$$\mu:\widetilde{\mathcal{P}} \to \widetilde{\mathcal{S}}\, , \ \ \ 
\mu_{\mathcal{A}} =  \widetilde{\mathcal{P}} \to {\mathcal{A}}_X\, $$
be the total space of the principal $G$-bundle over $\widetilde{\mathcal{S}}$. Then, 
$$S = \mathcal{S} \times_{\mathcal{A}_X} X\,, \ \  \overline{S} = \overline{\mathcal{S}} \times_{\mathcal{A}_X} X\, , \ \  \widetilde{S} = \widetilde{\mathcal{S}} \times_{\mathcal{A}_X} X\, , \ \  \widetilde{P} =  \widetilde{\mathcal{P}} \times_{\mathcal{A}_X} X\, .$$ Furthermore, since the map $\alpha$ is smooth, 
we find that $N_{\widetilde{S}/X}$ is the pullback of $N_{\widetilde{\mathcal{S}}/\mathcal{A}_X}$, and the splitting of $N_{\widetilde{S}/X}$ on $\widetilde{P}$ into line bundles is pulled back from the splitting of $N_{\widetilde{\mathcal{S}}/\mathcal{A}_X}$ on $\widetilde{\mathcal{P}}$. In other words, we have a Cartesian diagram:
\[
\begin{tikzcd}
\widetilde{P} \ar[r,"\alpha_P"] \ar[d,"\mu_X"] & \widetilde{\mathcal{P}} \ar[d,"\mu_{\mathcal{A}}"] \\ X  \ar[r,"\alpha"] & \ \mathcal{A}_X\, .
\end{tikzcd}
\]
Normally decorated strata classes on $\overline{S}$ have the form $\mu_{X*}\alpha_P^*(\gamma)$ for $\gamma \in \mathsf{CH}^\star(\widetilde{\mathcal{P}})$. As $\alpha,\alpha_P$ are smooth, $\mu_{X*}\alpha_P^*= \alpha^*\mu_{\mathcal{A}*}$. Therefore,
$$
R(X,D) \subset \alpha^*\mathsf{CH}^\star(\mathcal{A}_X)\, .
$$
In fact, the argument shows more precisely that 
$$
R(X,D) = \alpha^*R(\mathcal{A}_X,\mathcal{D})
$$
for $\mathcal{D} = \alpha(D)$ the corresponding divisor in $\mathcal{A}_X$. In other words, the logarithmic tautological ring of $(X,D)$, which is generated by the Chern roots of the normal bundles on the various monodromy torsors of the strata of $(X,D)$, is the pullback of the logarithmic tautological ring of $\mathcal{A}_X$, generated by the analogous constructions over the strata of $\mathcal{A}_X$. Thus, it suffices to show that the normally decorated strata classes of $\mathcal{A}_X$ generate the Chow ring of $\mathcal{A}_X$. We may thus reduce to proving the theorem for $\mathcal{A}_X$. 

So let $\gamma$ be a class in $\mathsf{CH}^\star(\mathcal{A}_X)$. We must show 
$$ \gamma \in R(\mathcal{A}_X,\mathcal{D})\, .$$
We may assume that $\gamma$ is supported on $\overline{\mathcal{S}}$ for 
some stratum $\mathcal{S}\subset \mathcal{A}_X$. Suppose, by induction, we have shown that every such class supported on a stratum $\overline{\mathcal{S}'}$ with $$\dim \mathcal{S}' < \dim \mathcal{S}$$ is in $R(\mathcal{A}_X,\mathcal{D})$. 
Suppose further that we can find a class $\delta \in R(\mathcal{A}_X,\mathcal{D})$ such that 
$\gamma$ equals  $\delta$ on $\mathcal{S}$. Then, 
$$\gamma-\delta\in  \mathsf{CH}^\star(\mathcal{A}_X)$$ is supported on lower dimensional strata and 
therefore lies in $R(\mathcal{A}_X,\mathcal{D})$, so that we have
$\gamma \in R(\mathcal{A}_X,\mathcal{D})$ as well. Thus, the induction hypothesis ensures that, for a given dimension $\dim \mathcal{S}$, we can remove strata $\mathcal{S}'$ with $\dim \mathcal{S}' < \dim \mathcal{S}$, and thus, it suffices to prove the statement with the additional assumption that $\mathcal{S}$ is closed in $\mathcal{A}_X$. Note that this reduction also suffices to handle the base of the induction: the minimal dimensional strata of $\mathcal{A}_X$ are automatically closed. 

Suppose then that $\gamma$ is a class supported on $\mathcal{S}$, and $\mathcal{S}$ is closed of codimension $n$ in $\mathcal{A}_X$. Then $\mathcal{S} \cong B(\mathbb{G}_m^n \rtimes G)$, $\mathcal{A}_X$ is a quotient of $[\mathbb{A}^n/\mathbb{G}_m^n]$ by an \'etale equivalence relation in a neighborhood of $\mathcal{S}$, and the monodromy torsor $\mu: \tilde{\mathcal{P}} \to \mathcal{S}$ is isomorphic to $B\mathbb{G}_m^n$.

We can use this to describe the normal bundle $N_{\mathcal{S}/\mathcal{A}_X}$ on $\mathcal{S}$: the data of this vector bundle on $\mathcal{S}$ is equivalent to specifying the bundle $\mu^* N_{\mathcal{S}/\mathcal{A}_X}$ on $\tilde{\mathcal{P}}$ together with a $G$-action. It is given by  
$$
\mu^* N_{\mathcal{S}/\mathcal{A}_X} = \oplus_{i=1}^n N_i := \oplus_{i=1}^n \mathcal{O}(\mathcal{D}_i)|_{\tilde{\mathcal{P}}}
$$
where $\mathcal{D}_i$ is the $i$-th hyperplane divisor in $[\mathbb{A}^n/\mathbb{G}_m^n]$. 
The monodromy group $G$ acts by permuting the hyperplanes $\mathcal{D}_i$ cutting out $S$, and this action lifts to a corresponding action permuting the direct summands $N_i$ above.
In particular, while the pullback $\mu^* N_{\mathcal{S}/\mathcal{A}_X}$ is a direct sum, the individual direct summands are in general not invariant under the $G$-action, and thus $N_{\mathcal{S}/\mathcal{A}_X}$ is not actually split on $\mathcal{S}$.

Still, on $\tilde{\mathcal{P}}$ we have that the classes $x_i:=c_1(N_i)$ form a generating set for the algebra
$$
\mathsf{CH}(\tilde{\mathcal{P}}) \cong \mathbb{Q}[x_1,\cdots,x_n]\,.
$$

On the other hand, the map $\mu$ gives an isomorphism 

$$
\mu^*: \mathsf{CH}(\mathcal{S}) \cong \mathsf{CH}(\tilde{\mathcal{P}})^G
$$
with inverse $\frac{1}{|G|}\mu_*$, since we are working with rational Chow groups. Thus $\gamma$ is the image of $\frac{1}{|G|}\mu^*\gamma$ under $\mu_*$, which is a $G$-invariant polynomial in the $x_i$. This shows that $\gamma$ is a normally decorated strata class, completing the proof. \qed

\vspace{10pt}


Theorem \ref{tautological=pp2} immediately implies that $\mathsf{R}^\star(X,D)\subset \mathsf{CH}^\star(X)$
is closed under the intersection product (a claim which was
left to the reader in Section \ref{tautc}). On the other hand, it is not immediate to see which piecewise polynomial corresponds to which normally decorated strata class. The precise correspondence between piecewise polynomials and normally decorated strata classes  has now been  carried out in \cite[Section 6]{HMPPS}.

\subsection{Proofs of Theorems
\ref{plbk} and \ref{pfwd}}
  \label{craf2}

Fix a normal crossings pair $(X,D)$ with Artin fan $\mathcal{A}_X$ and map 
$$\alpha:X \to \mathcal{A}_X\, .$$ Consider an arbitrary smooth log modification $$f: \widetilde{X} \to X\, $$ necessarily of the form $(\widetilde{X},\widetilde{D})$ with an associated map 
$$\widetilde{\alpha}: \widetilde{X} \to \mathcal{A}_{\widetilde{X}}\, .$$ By definition, the log modification $\widetilde{X} \to X$ is pulled back to $X$ from a log modification $\widetilde{\mathcal{A}}_{X} \to \mathcal{A}_X$ of Artin fans, and we have a diagram
\[
\begin{tikzcd}
\widetilde{X} \ar[r,"\beta"] \ar[d,"f"] & \widetilde{\mathcal{A}}_{X} \ar[d,"g"] \ar[r,"c"] & \mathcal{A}_{\widetilde{X}} \\ X  \ar[r,"\alpha"] & \ \mathcal{A}_X\ &
\end{tikzcd}
\]
with the square being Cartesian, the map $c$ proper, DM-type, \'etale and bijective, and $\widetilde{\alpha} = c \circ \beta$. By Theorem \ref{tautological=pp2}, $$\mathsf{R}^\star(X,D)=\alpha^*\mathsf{CH}^\star(\mathcal{A}_X) \ \ \
\text{and} \ \ \ \mathsf{R}^\star(\widetilde{X},
\widetilde{D})= \widetilde{\alpha}^*\mathsf{CH}^\star(\mathcal{A}_{\widetilde{X}})\, .$$
Since the map $c$ is proper, DM-type, \'etale and bijective, it induces an isomorphism 
$$
c^*:\mathsf{CH}^{\star}(\mathcal{A}_{\widetilde{X}}) \to \mathsf{CH}^{\star}(\widetilde{\mathcal{A}}_X)
$$
between rational Chow groups, and thus we also have 
$$
\mathsf{R}^\star(\widetilde{X},
\widetilde{D})= \beta^*\mathsf{CH}^\star(\widetilde{\mathcal{A}}_{X})\
$$
As $f_*{\beta}^*(\widetilde{\delta}) = \alpha^*g_*(\widetilde{\delta})\, ,$
we have
$$f_*\mathsf{R}^\star(\widetilde{X},\widetilde{D}) = \mathsf{R}^\star(X,D)\, ,$$ 
where we conclude equality instead of inclusion since $g_*$ is surjective.{\footnote{
{
For an interpretation of the pushforward $g_*$ in terms of piecewise polynomials, we refer the reader to \cite[Section 2.3]{BrionPP} where the toric setting is studied. These ideas are
used in the calculations of \cite{HMPPS}.
}
}}
Similarly, since
$f^*\alpha^*(\delta) = \beta^*g^*(\delta)$, we have $f^*\mathsf{R}^\star(X,D) \subset \mathsf{R}^\star(\widetilde{X},\widetilde{D})$.  \qed

Combining theorem \ref{tautological=pp2} with the techniques used in the proof above also provides a second proof
of Theorem \ref{cc77} based on
the study of the Artin fan. The crucial observation is the following. Suppose $(X,D)$ is a normal crossings pair with $D$ ``as simple as possible": $D$ is normal crossings in the Zariski topology, and the non-empty intersections of the branches of $D$ are connected. Equivalently, this means that $C(X,D)$ is the cone over an abstract simplicial complex, i.e. can be piecewise linearly embedded into a vector space. Then, the ring of piecewise polynomials on $C(X,D)$ has a global description in terms of the Stanley-Reisner ring: 

$$
\mathsf{PP}(C(X,D)) = \mathbb{Q}[x_r]/N
$$
where the variables $x_r$ range over the rays of $C(X,D)$, and $N$ is the ideal of non-faces, i.e. generated by monomials $x_{i_1}\cdots x_{i_k}$ ranging over the collections $i_1,\cdots,i_k$ of rays which do not form a cone in $C(X,D)$. A fortiori, this presentation implies that $\mathsf{CH}(\mathcal{A}_X)$ is generated by divisors. 

While the piecewise polynomials of a general $(X,D)$ do not admit this simple description, the observation is relevant in our context because any sufficiently fine log blowup of $(X,D)$ has this form. For example, the double barycentric subdivision $(\widehat{X},\widehat{D})$ of $(X,D)$ always has this form. Applying barycentric subdivision once on an arbitrary $C(X,D)$ produces a cone complex with no self-intersection (and thus no monodromy), but where two cones possibly share the same set of rays (i.e. the intersection of a set of branches of the divisor is disconnected). Applying barycentric subdivision a second time separates such cones, ensuring that each cone is uniquely characterized by its set of rays, and thus produces a cone complex $C(\widehat{X},\widehat{D})$ which is the cone over a simplicial complex.

\begin{corollary} \label{Cor:divlogCH}
We have $\mathsf{R}^\star(X,D)\subset \mathsf{divlogCH}^\star(X,D)$. 
\end{corollary}

\noindent {\em Proof.}
Let $(\widehat{X},\widehat{D})$ be
the 
 log blow-up 
 corresponding to the double barycentric subdivision,
 \[
\begin{tikzcd}
{\widehat{X}} \ar[r,"\widehat{\alpha}"] \ar[d,"f"] & \widehat{\mathcal{A}}_{X} \ar[d,"g"] \\ X  \ar[r,"\alpha"] & \ \mathcal{A}_X\, .
\end{tikzcd}
\]
As above, $\widehat{\mathcal{A}}_X$ is the relative Artin fan of $\widehat{X} \to X$, and the Artin fan $\mathcal{A}_{\widehat{X}}$ has the same rational Chow ring as $\widehat{\mathcal{A}}_X$. Let $\gamma \in \mathsf{R}^\star(X,D)$.
 By Theorem \ref{tautological=pp2},
 $\gamma\in \alpha^*\mathsf{CH}^\star(\mathcal{A}_X)$  and therefore
 $$f^*(\gamma) \in 
 \widehat{\alpha}^*\mathsf{CH}^\star(\mathcal{A}_{\widehat{X}})\, .$$
 Since
 $\mathsf{CH}^\star(\mathcal{A}_{\widehat{X}})$ is generated by divisors, we have
 $f^*(\gamma) \in \mathsf{divCH}^\star(\widehat{X})$. \qed 

\vspace{10pt}

The proof of Theorem \ref{tautological=pp2} immediately yields
a finer statement:
$ \mathsf{R}^\star(X,D)$ lies in the subalgebra generated by logarithmic divisors of the log blow-up associated to the second barycentric subdivision
of the Artin fan of $(X,D)$. In fact, the subalgebra generated by logarithmic divisors of the log blow-up associated to any log blowup $(\tilde{X},\tilde{D})$ with $C(\tilde{X},\tilde{D})$ the cone over a simplicial complex. The double barycentric subdivision of any normal crossings pair $(X,D)$ is always a canonical such choice, but for any given example, a much more efficient choice $(\widetilde{X},\widetilde{D})$ may be available.

\section{Pixton's formula for \texorpdfstring{$\lambda_g\in \mathsf{CH}^\star(\overline{\mathcal{M}}_g)$}{lambda\_g in CH*(Mbar\_g)}}

\subsection{Strata}
Pixton's formula for the double ramification cycle
$\mathsf{DR}_{g,A}\in \mathsf{CH}^g(\overline{\mathcal{M}}_{g,n})$ is
expressed as a sum over strata of $(\overline{\mathcal{M}}_{g,n},\partial \overline{\mathcal{M}}_{g,n})$  indexed by the set $\mathsf{G}_{g,n}$
of stable graphs.
We present here Pixton's formula with an emphasis on the special case
 $$\mathsf{DR}_{g,\emptyset} = (-1)^g\lambda_g \in \mathsf{CH}^g(\overline{\mathcal{M}}_{g})\, .$$
We refer the reader
to \cite{JPPZ,P} for a more detailed discussion about double ramification
cycles, stable graphs, Pixton's formula, and the relation to classical Abel-Jacobi theory.

\subsection{Weightings} \label{wwww}
Let $A=(a_1,\ldots,a_n)\in \mathbb{Z}^n$ satisfy $\sum_{i=1}^n a_i=0$. Let 
$$\Gamma \in \mathsf{G}_{g,n}$$
be a stable
graph\footnote{Here and in Pixton's formula in Section \ref{pffff}, we follow the notation of \cite[Sections 0.3 and 0.4]{JPPZ}.
The factors of 2 are treated equivalently but
slightly differently in \cite{BHPSS,JPPZ2}.
}
of genus $g$ with $n$ legs.
A {\em weighting}
of $\Gamma$ is a function on the set of half-edges,
$$ w:\H(\Gamma) \rightarrow \Z,$$
which satisfies the following three properties:
\begin{enumerate}
\item[(i)] $\forall h_i\in \L(\Gamma)$, corresponding to
 the marking $i\in \{1,\ldots, n\}$,
$$w(h_i)=a_i\ ,$$
\item[(ii)] $\forall e \in \E(\Gamma)$, corresponding to two half-edges
$h,h' \in \H(\Gamma)$,
$$w(h)+w(h')=0\, ,$$
\item[(iii)] $\forall v\in \V(\Gamma)$,
$$\sum_{v(h)= v} w(h)=0\, ,$$ 
where the sum is taken over {\em all} $n(v)$ half-edges incident to $v$.
\end{enumerate}
In the case $A=\emptyset$, the set of half-edges $\H(\Gamma)$
has no legs ($n=0$).

Let $r$ be a positive integer.
A {\em weighting mod $r$}
of $\Gamma$ is a function,
$$ w:\H(\Gamma) \rightarrow \{0,\ldots, r-1\},$$
which satisfies exactly properties (i-iii) above, but
with the equalities replaced, in each case, by the condition of 
{\em congruence $mod$ $r$}.
 The set $\mathsf{W}_{\Gamma,r}$ of such weightings $w$ is finite, with cardinality $r^{h^1(\Gamma)}$.

\subsection{Formula for double ramification cycles}
\label{pffff}
Let $A=(a_1,\ldots,a_n)\in \mathbb{Z}^n$ satisfy $\sum_{i=1}^n a_i=0$.
Let $r$ be a positive
integer.
We denote by
$$\P_g^{d,r}(A)\in R^d(\oM_{g,n})$$
the degree $d$ component of the tautological class 
\begin{multline} \label{gg99g}
\hspace{-10pt}\sum_{\Gamma\in \mathsf{G}_{g,n}} 
\sum_{w\in \mathsf{W}_{\Gamma,r}}
\frac1{|\Aut(\Gamma)| }
\, 
\frac1{r^{h^1(\Gamma)}}
\;
\xi_{\Gamma*}\Bigg[
\prod_{i=1}^n \exp(a_i^2 \psi_{h_i}) \cdot 
\\ \hspace{+10pt}
\prod_{e=(h,h')\in \E(\Gamma)}
\frac{1-\exp(-w(h)w(h')(\psi_h+\psi_{h'}))}{\psi_h + \psi_{h'}} \Bigg]\, .
\end{multline} 
in $R^*(\oM_{g,n})$.


The following fundamental polynomiality property of $\P_g^{d,r}(A)$  has been proven by Pixton,
see \cite[Appendix]{JPPZ}.

\begin{proposition}[Pixton] For fixed $g$, $A$, and $d$, the \label{pply}
class
$$\P_g^{d,r}(A) \in R^d(\oM_{g,n})$$
is polynomial in $r$ (for all sufficiently large $r$).
\end{proposition}

We denote by $\P_g^d(A)$ the value at $r=0$ 
of the polynomial associated to $\P_g^{d,r}(A)$ by Proposition~\ref{pply}. In other words, $\P_g^d(A)$ is the {\em constant} term of the associated polynomial in $r$. 
Pixton's
 formula for double ramification
cycles is
$$\mathsf{DR}_{g,A} = 2^{-g}\, \P_g^g(A)\, \in \mathsf{CH}^g(\oM_{g,n})\, .$$

 \subsection{Examples in the $A=\emptyset$ case}
  For the reader's convenience, we present here
  the first few examples{\footnote{The graphics are by  F. Janda.}}
  of Pixton's formula for $\lambda_g$ obtained 
  by calculating $(-1)^g \mathsf{DR}_{g,\emptyset}$.

Each labeled graph $\Gamma$ describes a moduli space $\oM_\Gamma$ (a product of moduli spaces associated with the 
vertices of $\Gamma$), a tautological class $\alpha \in R^*(\oM_\Gamma)$, and a natural map $$\xi_\Gamma: \oM_\Gamma \to 
\oM_g\, .$$
Our convention in the formulas below is that the graph $\Gamma$ represents the cycle class $(\xi_\Gamma)_* \alpha$. 
For instance, assume the graph carries no $\psi$-classes and the class $\alpha$ equals~1. Since the map $\xi_\Gamma$ is of degree $|\Aut (\Gamma)|$ onto its image, the cycle class represented by $\Gamma$ is then $|\Aut (\Gamma)|$ times the class of the image of~$\xi_\Gamma$.

\tikz{\coordinate (A) at (0,0); \coordinate (B) at (1,0); \coordinate (C) at (0.6,0.5);}

\tikzset{baseline=0, label distance=-3mm}
\def\NC{\draw (0,0.25) circle(0.25);}
\def\NL{\draw plot [smooth,tension=1.5] coordinates {(0,0) (-0.2,0.5) (-0.5,0.2) (0,0)};}
\def\NR{\draw plot [smooth,tension=1.5] coordinates {(0,0) (0.2,0.5) (0.5,0.2) (0,0)};}
\def\NN{\NL\NR}
\def\NNN{\NN \begin{scope}[rotate=180] \NR \end{scope}}
\def\NNNN{\NN \begin{scope}[rotate=180] \NN \end{scope}}
\def\NRS{\begin{scope}[shift={(B)}] \NR \end{scope}}
\def\NRD{\begin{scope}[rotate around={-90:(B)}] \NRS \end{scope}}
\def\DE{\draw plot [smooth,tension=1] coordinates {(0,0) (0.5,0.15) (1,0)}; \draw plot [smooth,tension=1.5] coordinates {(0,0) (0.5,-0.15) (1,0)};}
\def\DES{\begin{scope}[shift={(B)}] \DE \end{scope}}
\def\TE{\DE \draw (A) -- (B);}
\def\QE{\DE \draw plot [smooth,tension=1] coordinates {(A) (0.5,0.05) (B)}; \draw plot [smooth,tension=1.5] coordinates {(A) (0.5,-0.05) (B)};}
\def\T{\draw (0.2,0) -- (C) -- (B) -- (0.2,0);}
\def\TT{\draw (C) -- (B) -- (0.2,0); \draw plot [smooth,tension=1] coordinates {(0.2,0) (0.3,0.3) (C)}; \draw plot [smooth,tension=1] coordinates {(0.2,0) (0.5,0.2) (C)};}
\newcommand{\nn}[3]{\draw (#1)++(#2:3mm) node[fill=white,fill opacity=.85,inner sep=0mm,text=black,text opacity=1] {$\substack{\psi^#3}$};}
\renewcommand{\gg}[2]{\fill (#2) circle(1.3mm) node {\color{white}$\substack #1$};}

\paragraph{Genus~1.}
\begin{equation*}
 \lambda_1 = \frac 1{24} \tikz{\NC \gg{0}{A}}.
\end{equation*}

\paragraph{Genus~2.}
\begin{equation*}
\lambda_2 = 
  \frac 1{240} \tikz{\NC \nn{A}{130}{{}} \gg{1}{A}}
  + \frac 1{1152} \tikz{\NN \gg{0}{A}}.
\end{equation*}

\paragraph{Genus~3.}
\begin{align*}
\lambda_3 &= 
  \frac 1{2016} \tikz{\NC \nn{A}{130}{2} \gg{2}{A}}
  + \frac 1{2016} \tikz{\NC \nn{A}{130}{{}} \nn{A}{50}{{}} \gg{2}{A}}
  - \frac 1{672} \tikz{\DE \nn{A}{30}{{}} \gg{1}{A} \gg{1}{B}}
  + \frac 1{5760} \tikz{\NN \nn{A}{160}{{}} \gg{1}{A}} \\
  &
  - \frac{13}{30240} \tikz{\TE \gg{0}{A} \gg{1}{B}}
  - \frac 1{5760} \tikz{\NL \DE \gg{0}{A} \gg{1}{B}}
  + \frac 1{82944} \tikz{\NNN \gg{0}{A}}.
\end{align*}

\paragraph{Genus~4.}
\begin{align*}
\lambda_4 &= 
  \frac 1{11520} \tikz{\NC \nn{A}{130}{3} \gg{3}{A}}
  + \frac 1{3840} \tikz{\NC \nn{A}{130}{2} \nn{A}{50}{{}} \gg{3}{A}}
  - \frac 1{2880} \tikz{\DE \nn{A}{30}{2} \gg{1}{A} \gg{2}{B}}
  - \frac 1{3840} \tikz{\DE \nn{A}{30}{{}} \nn{A}{-30}{{}} \gg{1}{A} \gg{2}{B}}
  - \frac 1{1440} \tikz{\DE \nn{A}{30}{{}} \nn{B}{150}{{}} \gg{1}{A} \gg{2}{B}} \\
  &- \frac 1{1920} \tikz{\DE \nn{A}{30}{{}} \nn{B}{-150}{{}} \gg{1}{A} \gg{2}{B}}
  - \frac 1{2880} \tikz{\DE \nn{B}{150}{2} \gg{1}{A} \gg{2}{B}}
  - \frac 1{3840} \tikz{\DE \nn{B}{150}{{}} \nn{B}{-150}{{}} \gg{1}{A} \gg{2}{B}}
  + \frac 1{48384} \tikz{\NN \nn{A}{160}{2} \gg{2}{A}}
  + \frac 1{48384} \tikz{\NN \nn{A}{160}{{}} \nn{A}{110}{{}} \gg{2}{A}} \\
  &+ \frac 1{115200} \tikz{\NN \nn{A}{160}{{}} \nn{A}{20}{{}} \gg{2}{A}}
  + \frac 1{960} \tikz{\T \nn{B}{180}{{}} \gg{1}{B} \gg{1}{0.2,0} \gg{1}{C}}
  - \frac{23}{100800} \tikz{\TE \nn{A}{-30}{{}} \gg{2}{A} \gg{0}{B}}
  - \frac 1{57600} \tikz{\DE \NRS \nn{B}{20}{{}} \gg{2}{A} \gg{0}{B}} \\
&  - \frac 1{16128} \tikz{\DE \NRS \nn{B}{-150}{{}} \gg{2}{A} \gg{0}{B}}
  - \frac 1{16128} \tikz{\DE \NRS \nn{A}{-30}{{}} \gg{2}{A} \gg{0}{B}}
  - \frac 1{57600} \tikz{\DE \NRS \nn{B}{20}{{}} \gg{1}{A} \gg{1}{B}}
  - \frac 1{16128} \tikz{\DE \NRS \nn{B}{-150}{{}} \gg{1}{A} \gg{1}{B}} \\
  & - \frac 1{16128} \tikz{\DE \NRS \nn{A}{-30}{{}} \gg{1}{A} \gg{1}{B}}
  - \frac{23}{100800} \tikz{\TE \nn{A}{-30}{{}} \gg{1}{A} \gg{1}{B}}
  + \frac{23}{100800} \tikz{\TT \gg{2}{B} \gg{0}{0.2,0} \gg{0}{C}}
  + \frac{23}{50400} \tikz{\TT \gg{1}{B} \gg{1}{0.2,0} \gg{0}{C}}
  + \frac 1{16128} \tikz{\T \NRS \gg{0}{B} \gg{1}{0.2,0} \gg{1}{C}} \\
  &
  + \frac 1{115200} \tikz{\DE \DES \gg{1}{A} \gg{0}{B} \gg{1}{2,0}}
  + \frac 1{276480} \tikz{\NNN \nn{A}{20}{{}} \gg{1}{A}}
  - \frac{13}{725760} \tikz{\NL \TE \gg{1}{A} \gg{0}{B}}
  - \frac 1{138240} \tikz{\NL \NRS \DE \gg{1}{A} \gg{0}{B}} \\
  &
  - \frac{43}{1612800} \tikz{\QE \gg{1}{A} \gg{0}{B}}
  - \frac{13}{725760} \tikz{\NRS \TE \gg{1}{A} \gg{0}{B}}
  - \frac 1{276480} \tikz{\NRS \NRD \DE \gg{1}{A} \gg{0}{B}}
  + \frac 1{7962624} \tikz{\NNNN \gg{0}{A}}
\end{align*}

\subsection{Proof of Theorem \ref{xxx}}
\label{kkss33}
We analyze Pixton's formula in the $A=\emptyset$ case,
$$\lambda_g =
(-1)^g\mathsf{DR}_{g,\emptyset}   \in \mathsf{CH}^g(\overline{\mathcal{M}}_{g})
\, .$$
Since $A=\emptyset$, the sum \eqref{gg99g} is over stable graphs
$\Gamma\in \mathsf{G}_{g}$ corresponding to strata of 
$(\overline{\mathcal{M}}_{g},\partial \overline{\mathcal{M}}_{g})$.

\vspace{8pt}
\noindent $\bullet$ By the definition
of a {\em weighting mod $r$},
the weights 
$$w(h)\, ,\ w(h')$$  on the two halves of {\em every} separating 
edge $e$ of $\Gamma$ must both be $0$. The factor in Pixton's formula for 
$e$,
$$\frac{1-\exp(-w(h)w(h')(\psi_h+\psi_{h'}))}{\psi_h + \psi_{h'}}\, ,$$
then vanishes and kills the contribution of $\Gamma$ to $\P_g^g(\emptyset)$.
Therefore, nonvanishing terms in the sum \eqref{gg99g} must correspond
to graphs with {\em no} separating edges.

\vspace{8pt}
\noindent $\bullet$ Since $A=\emptyset$, the term
$$\prod_{i=1}^n \exp(a_i^2 \psi_{h_i})$$
drops out of  \eqref{gg99g}.

\vspace{8pt}
\noindent $\bullet$
The classes which do appear in \eqref{gg99g} are
the normal bundle terms $\psi_h+\psi_{h'}$ at
each edge of $\Gamma$.

\vspace{8pt}
Since the formula \eqref{gg99g} respects the automorphisms of the
stable graph $\Gamma$, we obtain the following result.

\begin{proposition} \label{pixL}
  The class $\lambda_g\in  \mathsf{CH}^g(\overline{\mathcal{M}}_{g})$
  is a sum of normally decorated classes associated to strata of
$(\overline{\mathcal{M}}_{g},\partial \overline{\mathcal{M}}_{g})$
  corresponding
  to stable graphs $\Gamma\in \mathsf{G}_{g}$ with no separating edges.
  \end{proposition}

  Theorem \ref{xxx} is then an immediate consequence of Proposition \ref{pixL}
  and Theorem \ref{cc77}. Proposition \ref{pixL} reflects a very special
  property of $\lambda_g$ obtained from Pixton's formula. 
 \qed
 
\vspace{10pt}

Since every
edge of every stable 
graph $\Gamma \in \mathsf{G}_g$
which appears
in Pixton's formula for $\lambda_g$
is non-separating, we actually have
$$\lambda_g \in 
\mathsf{R}^\star(\overline{\mathcal{M}}_g,\Delta_0)\, .$$
Theorem \ref{cc77} then implies
 a refinement of Theorem \ref{xxx},
$$\lambda_g \in 
\mathsf{divlogCH}^\star(\overline{\mathcal{M}}_g,\Delta_0)\, .$$
By applying Pixton's formula for
the double ramification cycle
$$\mathsf{DR}_{g,(0,\ldots,0)} = (-1)^g\lambda_g \in \mathsf{CH}^g(\overline{\mathcal{M}}_{g,n})\, ,$$
an identical argument yields 
  $$\lambda_g \in \mathsf{divlogCH}^\star(\overline{\mathcal{M}}_{g,n},\Delta_0)$$
  for $2g-2+n>0$.

\subsection{More general $\mathsf{DR}$ cycles}
\label{mgdrc}
  Let $A=(a_1,\ldots,a_n)$  be a vector of
  integers satisfying
  $\sum_{i=1}^n a_i =0$.
  Pixton's formula for  the double
 ramification cycle
 $$\mathsf{DR}_{g,A}\in \mathsf{R}^\star(\overline{\mathcal{M}}_{g,n})$$
 together with Theorem \ref{cc77} yields the
 following result (the proof is exactly the
 same as the proof of Theorem \ref{xxx}).
 \begin{theorem} \label{pss4}
 We have $\mathsf{DR}_{g,A}\in \underline{\mathsf{div}}\mathsf{logCH}^\star(\overline{\mathcal{M}}_{g,n})$
 where
$$ \underline{\mathsf{div}}\mathsf{logCH}^\star(\overline{\mathcal{M}}_{g,n}) \subset
\mathsf{logCH}^\star(\overline{\mathcal{M}}_{g,n})$$
is the subalgebra generated by logarithmic 
boundary divisors together with the cotangent
line classes $\psi_1,\ldots, \psi_n$.
 \end{theorem}

 Theorem \ref{pss4} provides half of the
 proof of Conjecture C concerning the
 lifted double ramification cycle 
 $\widetilde{\mathsf{DR}}_{g,A}$.
{There are now three proofs of the other half of the conjecture via three different approaches.
The first two
are by  Abel-Jacobi theory in \cite{HoSchwa} and by controlling the difference between $\mathsf{DR}_{g,A}$
and  $\widetilde{\mathsf{DR}}_{g,A}$ in an appropriate blowup of $\oM_{g,n}$ in \cite{MoRa}}. The third, presented in \cite{HMPPS}, proves the conjecture directly by giving a formula for (a representative of) $\widetilde{\mathsf{DR}}_{g,A}$ in terms of $\psi$-classes and piecewise polynomials.   

The special case $A=(0,\ldots,0)$ related to
the class
 $\lambda_g$ is simpler since
 no cotangent line classes appear at the markings
 in Pixton's formula.
Moreover, there is no
change in the lift for
$A=(0,\ldots,0)$:
$${\mathsf{DR}}_{g,(0,\ldots,0)}=
\widetilde{\mathsf{DR}}_{g,(0,\ldots,0)} \in 
 \underline{\mathsf{div}}\mathsf{logCH}^\star(\overline{\mathcal{M}}_{g,n})\, .$$

The $\omega^k$-twisted double ramification cycle
\cite{Hol} is also governed by Pixton's formula
\cite{BHPSS},
$$\mathsf{DR}^k_{g,A}\in \mathsf{R}^\star(\overline{\mathcal{M}}_{g,n})\, , \ \ \ \sum_{i=1}^n a_i = k(2g-2)\, .$$
The analogue of Theorem \ref{pss4}
can be proven for the $\omega^k$-twisted double ramification cycle, but the divisor subalgebra
of $\mathsf{logCH}^\star(\overline{\mathcal{M}}_{g,n})$
must  include $\kappa_1$ together with the
cotangent line classes
$\psi_i$ and the logarithmic
boundary divisors. Conjecture C can
then also be promoted
to a statement for the  
lifted
$\omega^k$-twisted double ramification cycle
(again including $\kappa_1$ in the subalgebra).

\subsection{Pixton's
generalized boundary strata classes}
In \cite{pixtongeneralboundary}, Pixton has defined a subalgebra of the tautological ring $\mathsf{R}^*(\overline{\mathcal{M}}_{g,n})$ spanned by \emph{generalized boundary strata classes}: tautological classes $[\Gamma]$ associated to prestable graphs $\Gamma$ of genus $g$ with $n$ legs. 

If $\Gamma$ is a semistable graph (every genus $0$ vertex is incident to at least two legs or half-edges),
then Pixton's definition
takes a simple form.
Let $\Gamma'$ be the stabilization of $\Gamma$.  The class $[\Gamma]$ is defined as a push-forward under the gluing map $\xi_{\Gamma'}$ of products of classes $\psi_1,\ldots, \psi_n$ and classes $\psi_h + \psi_{h'}$ for half-edges $(h,h')$ forming an edge of $\Gamma'$.
The analysis of Section \ref{kkss33}
then implies $$[\Gamma] \in  \underline{\mathsf{div}}\mathsf{logCH}^\star(\overline{\mathcal{M}}_g, \partial \overline{\mathcal{M}}_g)\,$$
in the semistable case.

Pixton's boundary class for more general unstable graphs
has $\kappa$ classes and will likely not lie in any version of
$\mathsf{divlogCH}^\star(\overline{\mathcal{M}}_g, \partial \overline{\mathcal{M}}_g)$.

\section{The bChow ring}\label{bch}

Let $X$ be a nonsingular variety. Given the additional data of a normal
crossings divisor $D \subset X$ we defined the log Chow ring of the pair
$(X,D)$. This is a variant of a much larger ring, the {\em bChow ring} of $X$. We define
\[
\mathsf{bCH}^\star(X) = \varinjlim_{Y \in \mathsf{B}(X)} \mathsf{CH}^\star(Y)\, ,
\]
where $\mathsf{B}(X)$ is the category of nonsingular blow-ups of $X$:
objects in $\mathsf{B}(X)$ are  proper birational maps
$$Y \rightarrow X $$
with $Y$ nonsingular
and morphisms in $\mathsf{B}(X)$ are proper birational maps over $X$.
For a longer introduction to the bChow ring, see \cite{HPS}.
Some of the ideas involved go back to 
papers of Shokurov \cite{Sh1,Sh2}. See
also Aluffi \cite{Al} for similar constructions.

Let $[Z\to X]$ and $[Y\to X]$ be objects of $\mathsf{B}(X)$.
If
$Z \to X$ factors as $$Z \to Y \to X\, ,$$
then there is a unique morphism from $[Z \to X]$ to $[Y \to X]$ in
$\mathsf{B}(X)$, and we call $Z \to X$ a {\em refinement}
of $Y \to X$. The transition maps in the above colimit are given by pullbacks
$$f^\star:\mathsf{CH}^*(Y) \to \mathsf{CH}^\star(Z)$$
for refinements $Z \stackrel{f}{\to} Y\to X$. 

Unlike, $\mathsf{logCH}^\star(X)$, the bChow ring does {\em not}
depend upon the choice of a normal crossings divisor $D\subset X$. However, given such a choice there is always a tower
of natural inclusions
$$\mathsf{CH}^\star(X) \subset \mathsf{logCH}^\star(X) \subset \mathsf{bCH}^\star(X)\,.$$
Since the centers of the blow-up are so restricted in the definition
of $\mathsf{logCH}^\star(X)$, we view $\mathsf{CH}^\star(X)$ and
$\mathsf{logCH}^\star(X)$ as relatively close in size. On the other hand,
$\mathsf{bCH}^\star(X)$ is very much larger.

Let $\mathsf{divbCH}^\star(X)$ be the subalgebra of $\mathsf{bCH}^\star(X)$
generated by divisors. More precisely,
\begin{equation*}
  \mathsf{divbCH}^\star(X) = \varinjlim_{Y \in \mathsf{B}(X)}
  \mathsf{divCH}^\star(Y)\, .
\end{equation*}
While the proof of the claim
$$\lambda_g \in \mathsf{divlogCH}^*(\overline{\mathcal{M}}_{g}
,\partial\overline{\mathcal{M}}_g)$$
depended upon special properties of $\lambda_g$,
the parallel bChow statement
$$\lambda_g \in \mathsf{divbCH}^*(\overline{\mathcal{M}}_{g})$$
immediately follows from a general result.

\begin{theorem}\label{rrr2}
For every nonsingular quasi-projective
variety{\footnote{The statement holds verbatim for nonsingular Deligne-Mumford stacks which
    admit finite resolutions of sheaves by vector bundles.}} $X$,
bChow 
is generated
by divisor classes,
$$\mathsf{divbCH}^*(X) = \mathsf{bCH}^\star(X)\, .$$
\end{theorem}

\noindent{\em Proof.} Let $\alpha \in \mathsf{CH}^\star(Y)$ for an object $[Y \to X]$ in $\mathsf{B}(X)$.
We will find a refinement $Z \to Y$ for which
$$f^*a \in \mathsf{divCH}(Z)\, .$$
Since $Y$ is nonsingular and quasi-projective, 
the Chern classes of vector bundles generate $\mathsf{CH}^\star(Y)$.
We can assume $\alpha = c_i(E)$ for a vector bundle $E$ on $Y$.
By \cite[Corollary 2]{H1}, there is a blow-up
$$g: W \to Y$$
where $W$ is nonsingular and
$g^*E$ contains a subline bundle $L$,
$$0 \rightarrow L \rightarrow g^*E \rightarrow g^*E/L \rightarrow 0 \, . $$
Applying the same argument to the quotient bundle
$g^*E/L$, we find inductively a nonsingular blow-up
$$f:Z \to Y$$ for
which $f^*E$ has a filtration with line bundles as quotients.
Therefore, $$f^*c_i(E) = c_i(f^*E)$$
is in $\mathsf{divCH}^\star(Z)$. 
\qed

\vspace{10pt}

The quasi-projective hypothesis is used only for vector bundle
resolutions. In fact, the hypothesis is not necessary. 
Theorem \ref{rrr2} can be proven locally near any cycle
$$ S \subset X\ $$
by successive blow-ups along nonsingular centers to 
resolve $S$ and appropriately modify the Chern classes
of the normal bundle of $S$. We leave the details for the interested reader.

\begin{appendix}

\section{The fourth cohomology group of \texorpdfstring{$\overline{\mathcal{M}}_{g}$}{Mbar\_g}} \label{cod2}
In the  proof of Theorem \ref{fff}, we require
 the equality{\footnote{We use,
as before, the complex grading on $\mathsf{RH}^\star$.}}
\begin{equation}\label{v445}
 \mathsf{H}^4(\overline{\mathcal{M}}_g) =
 \mathsf{RH}^2(\overline{\mathcal{M}}_g) 
 \,.
\end{equation}
for sufficiently large $g$,
In other words, 
the fourth cohomology group of $\overline{\mathcal{M}}_{g}$ is spanned by tautological classes for sufficiently
high $g$.

Equality \eqref{v445} 
was first proven by Edidin \cite{Edidin} 
for $g \geq 12$. Edidin bounded the Betti number  
$\mathsf{h}^4(\overline{\mathcal{M}}_g)$ from above and
then showed by intersection calculations that the span of 
the
tautological classes{\footnote{Edidin does not use
the language of tautological classes as we
now do, but all of his generators are
in fact tautological: they are given by the classes $\kappa_2, \kappa_1^2$, pushforwards of $\lambda$- and $\psi$-classes under boundary divisor gluing maps, and fundamental classes of strata of codimension $2$.}}  in codimension 2 achieves the required
rank. Edidin used the interior result
\begin{equation}
    \label{dd553}
    \mathsf{H}^4({\mathcal{M}}_g)
    =
\mathsf{RH}^2({\mathcal{M}}_g)  
\end{equation}
proven by Harer \cite{Harer} for $g\geq 12$. The interior
statement \eqref{dd553} was later proven for $g\geq 9$ by Ivanov \cite{Ivanov} and strengthened further to $g \geq 7$ by Boldsen \cite{Boldsen} which improved Edidin's bound. 

\begin{theorem}[\cite{Edidin}, \cite{Ivanov}, \cite{Boldsen}]
We have $\mathsf{H}^4(\overline{\mathcal{M}}_{g}) = \mathsf{RH}^2(\overline{\mathcal{M}}_{g})$ for $g \geq 7$.
\end{theorem}

\section{Computations in {\em admcycles}}
\subsection{Verification of Pixton's conjecture}
In \cite{Pixtonconj}, Pixton proposed a set of relations between tautological classes on the moduli spaces $\overline{\mathcal{M}}_{g,n}$ of stable curves. These were proven to hold in cohomology \cite{PPZ} and in Chow \cite{jandaEquivP1}. Pixton
furthermore conjectured that his relations span
the  \emph{complete} set of
relations among tautological classes. The relations were implemented by Pixton in the mathematical software SageMath \cite{SageMath} and later incorporated in the SageMath package {\em admcycles}. Assuming Pixton's conjecture, the software computes a basis of the $\mathbb{Q}$-vector spaces $\mathsf{R}^d(\overline{\mathcal{M}}_{g,n})$ and express 
tautological classes in the basis.

In Proposition \ref{Pixxx}, we state that Pixton's conjecture holds for the spaces
\[
\mathsf{R}^4(\overline{\mathcal{M}}_{4,1}) \text{ and }\mathsf{R}^5(\overline{\mathcal{M}}_{5,1})\,.
\]
Assuming the conjecture, {\em admcycles} computes the rank of these two spaces to be $191$ and $1371$ respectively. If the conjecture was false, the rank of one (or both) of the groups would have to be strictly smaller. However, using {\em admcycles}, we verify that the ranks of the intersection pairings
\[
\mathsf{R}^4(\overline{\mathcal{M}}_{4,1}) \otimes \mathsf{R}^{6}(\overline{\mathcal{M}}_{4,1}) \to \mathbb{Q} \text{ \ \ and\  \ }\mathsf{R}^5(\overline{\mathcal{M}}_{5,1}) \otimes \mathsf{R}^{8}(\overline{\mathcal{M}}_{5,1}) \to \mathbb{Q}
\]
are bounded from below by $191$ and $1371$ respectively. The rank
bounds are obtained
by taking generating sets of $\mathsf{R}^4(\overline{\mathcal{M}}_{4,1})$ and  $\mathsf{R}^5(\overline{\mathcal{M}}_{5,1})$ and computing the matrix of pairings with generators in $\mathsf{R}^6(\overline{\mathcal{M}}_{4,1})$ and $\mathsf{R}^8(\overline{\mathcal{M}}_{5,1})$ respectively.
For the rank bounds of pairing, we
do \emph{not} assume anything about the relations between the above generators, though we are allowed to use the known relations \cite{PPZ} to reduce the size of the generating sets. 

The computations were performed on a server
of the Max-Planck Institute for Mathematics in Bonn\footnote{The program ran on a single thread of the available CPU (Intel Xeon Prozessor E5-2667 v2) taking about 60 GB of RAM due to the large amounts of intermediate data to store (such as the list of Pixton's relations, sets of tautological generators, etc).}, taking two days in the case of $\overline{\mathcal{M}}_{4,1}$  and $31$ days for $\overline{\mathcal{M}}_{5,1}$. Without substantial improvements of the algorithm, it is thus unlikely that Pixton's conjecture can be verified in this way for significantly larger $g$, $n$, and $d$. We warmly thank the Max-Planck Institute for providing the computer infrastructure for our computations.

\subsection{Computations in proofs of Theorems \ref{rrr} and \ref{fff}}
Once we have verified Pixton's conjecture (as above{\footnote{Our verification method also then shows
$\mathsf{R}^d(\overline{\mathcal{M}}_{g,n}) =\mathsf{RH}^d(\overline{\mathcal{M}}_{g,n})$}}).
 for $\mathsf{RH}^d(\overline{\mathcal{M}}_{g,n})$, we can explicitly check whether 
\[\lambda_d \in \mathsf{RH}^d_{\leq k}(\overline{\mathcal{M}}_{g,n})\,.\]
Several such checks used in the proofs of Theorems \ref{rrr} and \ref{fff} were made using {\em admcycles}.

We provide below an example of the computation showing that the class $\lambda_3$ is not contained in the space 
$$\mathsf{divRH}^3(\overline{\mathcal{M}}_3) \subset
\mathsf{RH}^3(\overline{\mathcal{M}}_3)\,,$$
which is a 9-dimensional subspace of a 10-dimensional space.  We first create the list \texttt{divcl} of divisor classes on $\overline{\mathcal{M}}_3$, compute the set of triple products of such classes, and then take the span \texttt{divR} of the vectors representing them in a basis of $\mathsf{RH}^3(\overline{\mathcal{M}}_3)$. We verify that \texttt{divR} is $9$-dimensional inside the $10$-dimensional ambient space $\mathsf{RH}^3(\overline{\mathcal{M}}_3)$. Finally, we compute the class $\lambda_3$ and verify that the associated vector \texttt{Lv} is not contained in \texttt{divR}. 
\begin{verbatim}
sage: from admcycles import *
sage: divcl = tautgens(3,0,1)
sage: divp = [a*b*c for a in divcl for b in divcl for c in divcl]
sage: divR = span(u.toTautbasis() for u in divp)
sage: (divR.rank(), divR.degree())
(9, 10)
sage: L = lambdaclass(3,3,0)
sage: Lv = L.toTautbasis()
sage: Lv in divR
False
\end{verbatim}

\subsection{Proof of Proposition \ref{noT}} \label{app:lambda2calculation}
We record below the computation in {\em admcycles} used in the proof of Proposition \ref{noT}. We create the classes $\lambda_2, [\Delta_0], [B]$ and $[C]$ and represent the class defined by
\[2 \lambda_2 - x \cdot [\Delta_0]^2 - y \cdot [B] - z \cdot [C]\]
in the vector  \texttt{diff} with respect to a basis of $\mathsf{CH}^2(\overline{\mathcal{M}}_2) = {\mathsf{R}}^2(\overline{\mathcal{M}}_2)$. We then solve the equation \texttt{diff=0} to find the formula for $x$ and $y$ in terms of the variable $z$  used in the proof.

We remark that in the definition of the class \texttt{Delta0} we need to divide by $2$ since this is the degree of the gluing morphism parameterizing the boundary divisor $\Delta_0$.
\begin{verbatim}
sage: from admcycles import *
sage: lambda2 = lambdaclass(2,2,0)
sage: Delta0 = 1/2 * irrbdiv(2,0) 
sage: gammaB = StableGraph([0],[[1,2,3,4]],[(1,2),(3,4)])
sage: B = gammaB.boundary_pushforward()
sage: gammaC = StableGraph([0,1],[[1,2,3],[4]],[(1,2),(3,4)])
sage: C = gammaC.boundary_pushforward()
sage: x, y, z = var('x, y, z')
sage: diff = (2*lambda2 - x*Delta0^2 - y*B - z*C).toTautbasis()
sage: diff
(476*x + 1824*y - 96*z - 3, -144*x - 576*y + 24*z + 1)
sage: solve([diff[i]==0 for i in (0,1)], x,y,z)
[[x == r1 - 1/120, y == -5/24*r1 + 11/2880, z == r1]]
\end{verbatim}

\end{appendix}

\vspace{8pt}

\noindent Departement Mathematik, ETH Z\"urich\\
\noindent samouil.molcho@math.ethz.ch

\vspace{8pt}

\noindent Departement Mathematik, ETH Z\"urich\\
\noindent rahul@math.ethz.ch

\vspace{8pt}

\noindent Mathematisches Institut der Universit\"at Z\"urich\\
\noindent johannes.schmitt@math.uzh.ch

\end{document}